\documentstyle[amssymb,amsfonts,12pt,psfig]{amsart}

\newtheorem{theorem}{Theorem}[section]
\newtheorem{lemma}[theorem]{Lemma}
\newtheorem{corollary}[theorem]{Corollary}
\newtheorem{proposition}[theorem]{Proposition}

\def\C{{\mbox{\rm\kern.24em
\vrule width.03em height1.43ex depth-.052ex \kern-.26em C}}}
\def\QSet{\mbox{\rm\kern.24em
\vrule width.03em height1.48ex depth-.051ex \kern-.26em Q}}
\def\Z{{\bf Z}}
\def\R{{\mbox{\rm I\kern-.22em R}}}
\def\P{{\bf P}}

\def\Q{{\bf Q}}
\def\T{{\bf T}}

\def\size{{\rm size}}

\def\mean{{\rm mean}}

\def\I{{\bf I}}
\def\J{{\bf J}}
\def\eps{\varepsilon}

\def\Re{{\rm Re}}
\def\Im{{\rm Im}}
\def\S{{S}}
\def\CZ{{\rm T}}

\def\111{\gamma}

\def\Tree{{\rm Tree}}

\def\be#1{\begin{equation}\label{#1}}
\def\bas{\begin{align*}}
\def\eas{\end{align*}}
\def\bi{\begin{itemize}}
\def\ei{\end{itemize}}
\newenvironment{proof}{\noindent {\bf Proof} }{\endprf\par}
\def \endprf{\hfill  {\vrule height6pt width6pt depth0pt}\medskip}
\def\emph#1{{\it #1}}

\title{Carleson measures, trees, extrapolation, and $\CZ(b)$ theorems}

\author{P. Auscher}
\address{LAMFA, CNRS UMR 6140, Facult\'e de Math\'ematiques et d'Informatique ,
Universit\'e de Picardie-Jules Verne, 
80039 Amiens CEDEX 1, France 
}
\email{Pascal.Auscher@@u-picardie.fr}

\author{S. Hofmann}
\address{Mathematics Department, University of Missouri, Columbia MO, 65211 USA}
\email{hofmann@@math.missouri.edu}

\author{C. Muscalu}
\address{Department of Mathematics, UCLA, Los Angeles CA 90095-1555}
\email{camil@@math.ucla.edu}

\author{T. Tao}
\address{Department of Mathematics, UCLA, Los Angeles CA 90095-1555}
\email{tao@@math.ucla.edu}

\author{C. Thiele}
\address{Department of Mathematics, UCLA, Los Angeles CA 90095-1555}
\email{thiele@@math.ucla.edu}

\begin{document}

\begin{abstract}  The theory of Carleson measures, stopping time arguments, and atomic decompositions has been well-established in harmonic analysis.  More recent is the theory of phase space analysis from the point of view of wave packets on tiles, tree selection algorithms, and tree size estimates.  The purpose of this paper is to demonstrate that the two theories are in fact closely related, by taking existing results and reproving them in a unified setting.  In particular we give a dyadic version of extrapolation for Carleson measures, as well as a two-sided local dyadic $\CZ(b)$ theorem which generalizes earlier $\CZ(b)$ theorems of David, Journe, Semmes, and Christ.
\end{abstract}

\maketitle

\tableofcontents

\section{Introduction}\label{introduction-sec}

The purpose of this article is to demonstrate the close connection between two sets of techniques in harmonic analysis: the theory of Carleson measures and related objects, and the theory of trees and related objects. 

A Carleson measure is a positive measure $\mu$ on the upper half space such
   that $\mu(I \times (0,\ell(I)) \lesssim |I|$ for every cube $I \subseteq \R^n $ with side length $\ell(I)$.  There
   is also an analogous notion for domains more general than the half-space, as
   well as a discrete version:  if $\mu$ is a mapping from dyadic cubes into the
   non-negative reals, then $\mu$ satisfies a (discrete) Carleson measure
   condition $\sum_{I\subseteq J} \mu_I \lesssim |J|$ for every dyadic cube $J,$
   where the sum runs over all dyadic
   sub-cubes of $J.$  Carleson measures  are intimately connected with many
   aspects of harmonic analysis, including non-tangential behavior of functions
   in the half-space (or in a domain) (see e.g. \cite{stein}), $H^p$ theory and BMO
   \cite{fs}, boundedness of singular integrals, square functions and maximal
   functions (e.g. \cite{cj}, \cite{dj}, \cite{j1}, \cite{semmes}, \cite{k1}, \cite{k2}), geometric measure theory 
   (e.g. \cite{ds}, \cite{j1,j2}), and PDE (e.g. \cite{kato-extrap}, \cite{fkp}, \cite{hl},
   \cite{lm}).  Moreover, via their connection with the theory of trees, Carleson
   measures have played a significant role in recent work on Bilinear Singular
   Integrals \cite{laceyt1, laceyt2, cct, mtt:biest, thiele}, and (rather
appropriately!) Carleson's theorem on a.e.   convergence of Fourier Series \cite{fefferman}, \cite{laceyt-carleson}.  In these latter connections it is more convenient to work in the phase plane than in the Carleson half-space, and we have deliberately chosen our notation to reflect this fact.

This article is mainly expository.  Apart from one main new result (a local $\CZ(b)$ theorem), we shall mostly take existing results (atomic decompositions, paraproduct estimates, Carleson embedding) and re-prove them in a framework which unifies both the Carleson measure theory and the theory of trees and tiles.  (As such there is some overlap with the recent lecture notes in \cite{pereyra}).

Since this is an expository article, we shall simplify matters and only work in one dimension $\R$.  Also, we shall mostly work in the dyadic setting instead of the continuous one, to avoid issues such as rapidly decreasing tails or use of the Vitali covering lemma.  Thus, our results will be phrased using dyadic intervals and the Haar basis instead of arbitrary intervals and Gaussians (or similar smooth kernels).  However most of our results have continuous analogues (see e.g. \cite{gj} for a comparison between dyadic and continuous harmonic analysis).  We also will truncate all our spaces to be finite-dimensional to avoid technicalities.

The paper is organized as follows.  After setting up the notation of dyadic Carleson measures and BMO, Haar wavelets, and tiles and trees, we will give a quick review of the standard ``$L^\infty$'' theory of BMO (i.e. measuring the ways in which BMO is close to $L^\infty$), but from the perspective of trees and tiles. As part of this $L^\infty$ theory, we give a trees-based proof of the (dyadic analogue of the) extrapolation lemma for Carleson measures developed recently in \cite{kato-extrap}, \cite{hl}, \cite{lm}.  We also give an alternate proof of
the extrapolation lemma due to John Garnett.

We then show how BMO is also useful in ``$L^p$'' contexts, mainly through a BMO version of the Calder\'on-Zygmund decomposition.   This type of lemma is used often in the recent work on Carleson's theorem and the bilinear Hilbert transform, and is implicit in earlier work on Carleson measures and similar objects; we illustrate this by using the BMO Chebyshev inequality to re-prove the standard atomic decomposition of $H^p$.  

Next, we prove the Carleson embedding theorem and give its usual applications to paraproduct estimates and the $\CZ(1)$ theorem.  We also give a short proof of the boundedness of paraproducts below $L^1$; the proof is more direct than earlier proofs in that one does not go explicitly through the $\CZ(1)$ theorem.

Finally, we consider Calder\'on-Zygmund operators.  We prove a two-sided local $\CZ(b)$ theorem which generalizes the existing local and global $\CZ(b)$ theorems (\cite{djs}, \cite{at}, \cite{christ-tb}, \cite{semmes}); for instance, we can prove the standard global $\CZ(b)$ theorem assuming that $b$ is only in $BMO$ rather than $L^\infty$. 

The $\CZ(b)$ Theorem, in its various guises, has its roots in a question
posed by Yves Meyer, who asked whether the $\CZ(1)$ Theorem of David and
Journ\'e \cite{dj} (see also Chapter 6 below)
remains true if the constant function $1$ is replaced
by some function $b \in L^{\infty}$ with ${\rm Re} \, b \geq \delta $ (such
$b$ are said to be ``accretive").  The question was motivated by its
applicability to the $L^2$ boundedness of the Cauchy integral operator
on a Lipschitz graph.  Indeed, if $\Gamma$ denotes the graph, in the
plane, of a real-valued Lipschitz function $A,$ then by Cauchy's
theorem, we have that in the sense of BMO (that is, modulo constants),
$$ 0 = p.v. \, \int_{\Gamma} \frac{1}{z-w} dw,$$
for $z \in \Gamma.$
But in graph co-ordinates, this amounts to saying that (again in the
sense of BMO),
$$ 0 = p.v. \, \int_{-\infty}^{\infty} \frac{1 + iA'(y)}{x-y +
i(A(x)-A(y))} dy \equiv T(b) (x),$$
where $b$ is the accretive function $1 + iA' ,$ and $T$ is the singular
integral operator naturally associated to the antisymmetric
Calderon-Zygmund kernel $K(x,y) =\left (x-y + i(A(x)-A(y))\right)^{-1}.$

The $L^2$ boundedness of $\CZ$, and hence also that of the Cauchy integral
operator
$$ C_{\Gamma} f (x) \equiv  p.v. \, \int_{\Gamma} \frac {f(w)}{z-w} dw,
$$ thus follows from an analogue of the $\CZ(1)$ Theorem in which the
condition $\CZ(1)$,  $\CZ^*(1) \in BMO$ is replaced by the condition 
$\CZ(b) = 0 = \CZ^* (b),$ for some accretive function $b.$  Just such a result was
proved by McIntosh and Meyer \cite{McM}, who consequently obtained an
alternative proof of their earlier joint result with Coifman \cite{CMcM}
concerning the $L^2$ boundedness of the operator $C_{\Gamma}.$

The ``$\CZ(b)$ Theorem" of \cite{McM} was generalized by David, Journe and
Semmes \cite{djs} to allow
$\CZ(b), \CZ^*(b) \in BMO$ (indeed, they allowed other generalizations as
well, for example that there could be
two different accretive functions $b_1,b_2$ such that $\CZ(b_1), \CZ^*(b_2)
\in BMO,$ and moreover that the pointwise accretivity condition could be
relaxed to a condition holding on various sorts of averages - see, e.g., 
the notion
of ``pseudo-accretivity" defined in Section 6.1 below).  

This led to a proof of the $\CZ(b)$ Theorem  by constructing Haar 
wavelets adapted to the function $b$ \cite{cjs} (we shall base our 
proof on a variation of these adapted Haar wavelets).

A very simple
proof of a ``one-sided version" of the $\CZ(b)$ Theorem was obtained by
Semmes \cite{semmes}, who observed that in the special case $\CZ(b) \in BMO, \,
\CZ^*(1) = 0,$ one can readily show that $\CZ(1) \in BMO,$ thus reducing
matters to the $\CZ(1)$ Theorem.
It is worth noting that a suitable adaptation of Semmes's argument is
applicable to the solution of the square root problem of Kato.  Indeed,
one of the present authors (Auscher), along with Tchamitchian
\cite{AT},  formulated a version of the $\CZ(b)$ Theorem whose proof was based
upon the argument of \cite{semmes}, and which was
subsequently used to
solve the Kato problem in higher dimensions \cite{kato}, \cite{HLMc},
\cite{kato-2}.
We further note that there are local versions of the $\CZ(b)$ Theorem, 
due to M. Christ \cite{christ-tb} (cf.
Theorem \ref{ctb} below), 
which also have interesting
applications, namely to questions of analytic capacity; see for instance
\cite{verdera} for further discussion. 

This work was conducted at U. Missouri, UCLA, and the Centre for Mathematics and its Applications (CMA) at ANU.  The authors are particularly grateful to CMA for their warm hospitality during the visit of three of us (PA, CT, TT).  SH is supported by NSF grant DMS 0088920.  CM is supported by NSF grant DMS 0100796.
TT is a Clay Prize Fellow and is supported by a grant from the Packard Foundation.  CT is supported by a Sloan fellowship and NSF grants DMS 9985572 and DMS 9970469.  The authors are indebted to Stephanie Molnar, John Garnett, Joan Verdera and the referee for many helpful corrections and comments.

\section{Notation}\label{notation-sec}

We use $A \lesssim B$ to denote the estimate $A \leq CB$ for some absolute constant $C$ which may vary from line to line.  

If $E$ is a set, we use $|E|$ to denote the Lebesgue measure of $E$.  We will always be ignoring sets of measure zero, thus we only consider two sets $E$, $F$ to be intersecting if $|E \cap F| > 0$. 

Although our functions may be complex valued, we shall use the real inner product
$$ \langle f, g \rangle := \int f(x) g(x)\ dx$$
throughout.

\subsection{Tiles and trees}

We shall be working with dyadic intervals throughout the paper.  The number of dyadic intervals is infinite, but to simplify the arguments we shall restrict ourselves to a finite set on the half-line; in applications, this restriction can always be removed by a standard translation and limiting argument.  Specifically, we fix a large integer $M > 0$; none of our estimates will depend on $M$.  We define \emph{dyadic interval} to an interval\footnote{We will be careless about whether our intervals are closed, half-open, or open because of our convention of ignoring sets of measure zero.} of the form $I = [j 2^k, (j+1) 2^k]$, where $j, k$ are integers such that $-M \leq k \leq M$ and $I \subseteq [0, 2^M]$.  Let $\I$ denote the set of all dyadic intervals; observe that $\I$ is finite.  All sums and unions involving $I$ or $J$ will be assumed to be over $\I$ unless
otherwise specified.  If $f$ is a function on $\R$, we define $[f]_I := \frac{1}{|I|} \int_I f$ to denote the mean of $f$ on $I$.  We use $2I$ to denote the parent\footnote{On the non-dyadic theory $2I$ is often used to denote the interval with the same center as $I$ but twice the length; this can be thought of as a non-dyadic version of the parent of $I$.  However, in this paper we use $2I$ to exclusively refer to the dyadic parent of $I$, i.e. the unique dyadic interval of twice the length which contains $I$.} of $I$, and $I_l$, $I_r$ to denote the left and right children of $I$ (these are undefined if $|I|=2^M$ or $|I|=2^{-M}$ respectively).  We refer to the intervals $I_l$ and $I_r$ as \emph{siblings}.

Since our dyadic intervals have been restricted to a finite set, all norms will automatically be finite and all stopping time processes will automatically terminate.  This allows us to avoid some minor technicalities in our arguments, although it also means that we occasionally have to treat the smallest scale $|I| = 2^{-M}$ or the largest scale $|I| = 2^M$ a little differently from all the other scales.

A major advantage of the dyadic setting is the \emph{nesting property}: if $I, J$ are dyadic intervals which intersect each other, then either $I \subseteq J$ or $J \subseteq I$.  In particular, for any collection of dyadic intervals, the maximal intervals in this collection will always be disjoint.

The theory of Carleson measures is usually set in the \emph{upper half-space}
$$ \R^2_+ := \{ (x,t): x \in \R, t \in \R^+ \}.$$
Actually, because of our truncation parameter $M$ we will work in the compact subset
$$ (\R^2_+)_M := \{ (x,t): x \in [0, 2^M], t \in [2^{-M}, 2^{M-1}] \}.$$
The variable $x$ represents spatial position, while the $t$ variable represents time, wavelength, or spatial scale.  For every dyadic interval $I \in \I$, we let $l(I) = |I|$ denote the side-length of $I$, and define the \emph{Carleson box} $Q(I) \subset (\R^2_+)_M$ by
$$ Q(I) := I \times [2^{-M}, l(I)]$$
and the \emph{Whitney box} $Q^+(I) \subset Q(I)$ by
$$ Q^+(I) := I \times [\frac{l(I)}{2}, l(I)].$$
We remark that we have the partition 
$Q(I) = \bigcup_{J: J \subseteq I} Q^+(J)$.

Meanwhile, the theory of trees and tiles is usually set in \emph{phase space}
$$ \R^2 := \{ (x,\xi): x \in \R, \xi \in \R \}.$$
Because of our truncation, and because we are in the dyadic setting, we will instead work in the region\footnote{In truth, we are working not with the Euclidean field $\R$, but with the Walsh field $\R^+ \equiv (\Z_2)^\Z$.  See e.g. \cite{thiele-walsh}.}
$$ (\R^2)_M := \{ (x,\xi): x \in [0, 2^M], \xi \in [0, 2^M] \}.$$
The variable $x$ represents spatial position, while $\xi$ represents frequency.
A \emph{Heisenberg tile} or simply \emph{tile} is a rectangle in $\R^2$ of the form $P := I_P \times \omega_P$, where $I_P$ and $\omega_P$ are dyadic intervals such that $|P| = |I_P| |\omega_P| = 1$.

If $P$ and $Q$ are tiles, we say that $P \leq Q$ if $P$ intersects $Q$ and $I_P \subseteq I_Q$.  This is a partial order on tiles.

\begin{figure}[htbp]\label{tiles-fig} \centering
\ \psfig{figure=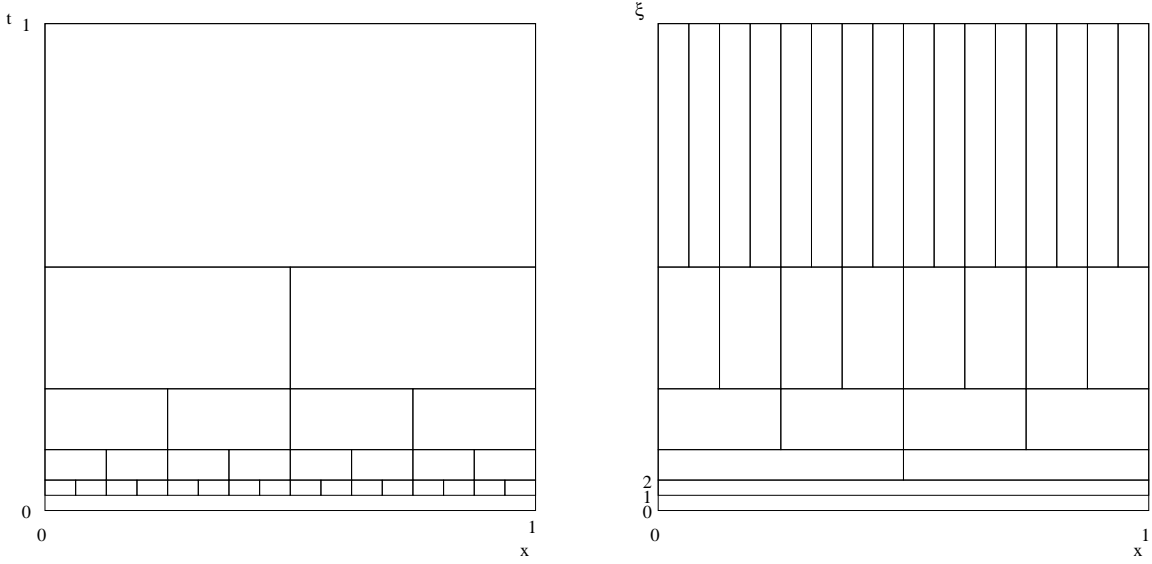,width=6in}
\caption{The geometry of the Carleson half-plane (partitioned into Whitney boxes) and phase space (partitioned into non-lacunary tiles).  The heuristic $t = 1/\xi$ provides a one-to-one correspondence between the two partitions. }
\end{figure}

If $I$ is a dyadic interval, we define the \emph{lacunary tile} $P^+(I)$ by
$$ P^+(I) := I \times [\frac{1}{l(I)}, \frac{2}{l(I)}]$$
and the \emph{non-lacunary tile} $P^0(I)$ by
$$ P^0(I) := I \times [0, \frac{1}{l(I)}].$$
Let $\P^+$ denote the set of all lacunary tiles, and $\P^0$ the set of all non-lacunary tiles. We define $P^+(I) \leq' P^+(J)$ if and only if $P^0(I) \leq P^0(J)$, or equivalently if $I \subseteq J$.  Thus $\leq'$ is a partial ordering on $\P^+$. Of course, there are many tiles which are not in either of these two sets, and many results in this paper can be extended to general tiles.  However, for simplicity we shall mostly restrict ourselves to the lacunary and non-lacunary tiles.  We write $[f]_P$ as shorthand for the averages $[f]_{I_P}$.

If $P^+(I)$ is a lacunary tile, we define the \emph{parent} $2P^+(I)$ of $P^+(I)$ by $2P^+(I) := P^+(2I)$.  Similarly define $2P^0(I) := P^0(2I)$.

A \emph{lacunary\footnote{Non-lacunary trees $T \subseteq \P^0$ are also useful in the study of the bilinear Hilbert transform and Carleson's operator; more precisely, when treating the bilinear Hilbert transform $\langle B(f,g), h \rangle$ one uses a triple of trees associated to $f$, $g$, $h$ respectively, with two of the trees lacunary and the third non-lacunary (but possibly with a non-zero frequency origin).  Similarly when treating the Carleson operator $\langle {\cal C}_{N(x)} f, \chi_E \rangle$ one uses a pair of trees associated to $f$ and $\chi_E$ respectively, with one lacunary and one non-lacunary.  See \cite{laceyt2}, \cite{laceyt-carleson}.  However we will not use non-lacunary trees explicitly in this paper, although they appear implicitly in Lemma \ref{mean-select} and in the paraproduct theory.} tree} (henceforth abbreviated as \emph{tree}) is a collection $T \subseteq \P^+$ of lacunary tiles with a \emph{top tile} $P_T \in T$, such that $P \leq' P_T$ for all $P \in T$.  We use $I_T$ as short-hand for $I_{P_T}$.  If $P \in \P^+$, we define the \emph{complete tree} $\Tree(P)$ to be the tree
$$ \Tree(P) := \{ Q \in \P^+: Q \leq' P \}$$
with top $P$.  We sometimes write $\Tree(I)$ for $\Tree(P^+(I))$.  Note that every tree $T$ lies inside a complete tree $\Tree(P_T)$.  If $T$ is a tree inside a collection $\P$ of tiles, we say that $T$ is \emph{complete with respect to $\P$} if $T = \Tree(P_T) \cap \P$.

Let $\alpha > 0$ and $T$ be a tree.  We define an \emph{$\alpha$-packing} of $T$ to be a set $\P \subset T$ of tiles such that
$$ \sum_{P \in \P} |I_P| \leq \alpha |I_T|.$$
We say that $\P$ is a \emph{uniform $\alpha$-packing}\footnote{This is roughly equivalent to $\sum_{P \in \P} \chi_{I_P}$ having a BMO norm bounded by $\alpha$.} of $T$ if
$$ \sum_{P \in \P: I_P \subset J} |I_P| \leq \alpha |J|$$
for all dyadic intervals $J$.

If $\alpha < 1/2$ and $\P$ is an $\alpha$-packing of $T$, observe that the parent tiles $2\P := \{ 2P: P \in \P\}$ form a $2\alpha$-packing of $T$.  Similarly if $\P$ is a uniform $\alpha$-packing of $T$, then $2\P$ is a uniform $2\alpha$-packing of $T$.

We say that a collection of lacunary tiles $\P$ is \emph{convex} if for every pair of tiles $P_1 \leq' P_2$ in $\P$, the set $\{ P \in \P^+: P_1 \leq' P \leq' P_2 \}$ is also contained in $\P$.  We will usually be dealing with convex trees in this paper.

The correspondence between the upper half-space and phase space is given by the heuristic formula\footnote{To be completely precise, one would have to adjust this formula when $\xi \leq 2^{1-M}$, but as this is only a heuristic anyway we will not bother to do this.}
\be{correlation}
t = 1/\xi;
\end{equation}
in other words, frequency is the reciprocal of wavelength.  This correspondence identifies Whitney boxes $Q^+(I)$ with lacunary tiles $P^+(I)$, and identifies a Carleson box $Q(I)$ with the complete tree $\Tree(I)$.  (Incomplete trees $T$ are identified with the portion of a Carleson box above a ``dyadic Lipschitz graph'', cf. \cite{kato-extrap}).  Note how this correspondence clearly gives a privileged position to the frequency origin $\xi = 0$.

The thesis of this paper is that the theory of Carleson measures can be equated with the theory of lacunary tiles.  The theory of general tiles - which is needed for applications such as Carleson's theorem and the bilinear Hilbert transform, in which the frequency origin plays no distinguished role - can then be thought of as a generalization of Carleson measure theory\footnote{In our paper, we will only need tiles which are centered at or near the frequency origin, in which case it does not particularly matter whether we use the Carleson half-plane or the phase plane.  However, we have chosen to use phase space notation (using frequency $\xi$ instead of wavelength $t$) as this is more compatible with the more general theory of multilinear operators such as the bilinear Hilbert transform (or the Carleson maximal operator), which are invariant under translations of the frequency variable.  Note that the modulation operation $f \mapsto e^{2\pi i\xi_0 x} f$ can be represented easily in the phase plane as a translation by $\xi_0$ in the $\xi$ variable, but is not so elegantly representable in the Carleson half-plane.  Nevertheless, we will not need to modulate in frequency in this paper, so the Carleson viewpoint and the phase space viewpoint are essentially equivalent here.}.

\begin{figure}\label{fig1}
\begin{tabular}{|l|l|} \hline
Tiles and trees & Carleson measures \\\hline
Phase space $\R^2$ & Upper half-plane $\R^2_+$ \\
Lacunary tile $P^+(I)$ & Whitney box $Q^+(I)$ \\
Non-lacunary tile $P^0(I)$ & ``Tower'' $I \times [l(I),\infty)$ \\
Complete tree $\Tree(I)$ & Carleson box $Q(I)$ \\
Convex tree & Carleson box above a Lipschitz graph \\
Size $\| \mu \|_{\size(T)}$ & Normalized mass $\mu(Q(I))/|I|$ \\
Bounded maximal size & Carleson measure (or BMO function)\\
$\xi$ & $1/t$ \\
\hline
\end{tabular}
\caption{A partial dictionary between tree terminology, and Carleson measure terminology.  In our paper the two viewpoints are essentially equivalent, however the phase space viewpoint is better adapted to handle more general situations where one needs to modulate in frequency.  Conversely, Carleson measures are 
better adapted to complex analysis applications.}
\end{figure}

\subsection{Size and Carleson measures}

Let $T$ be a convex lacunary tree, and suppose that we have a function $a: T \mapsto \R^+$ assigning non-negative numbers to each tile in $T$. We define the \emph{size} of $a$ on $T$ by
\be{size-def}
\| a \|_{\size(T)} := \frac{1}{|I_T|} \sum_{P \in \T} a(P).
\end{equation}
Now let $\P$ be any collection of lacunary tiles, and let $a: \P \mapsto \R^+$.  We define the \emph{maximal size} of $a$ on $T$ by
\be{size-max-def}
\| a \|_{\size^*(\P)} := \sup_{T \subset \P} \| a \|_{\size(T)}
\end{equation}
where $T$ ranges over all convex lacunary trees in $\P$; we adopt the convention that $\|a\|_{\size^*(\P)} = 0$ if $\P$ is empty.   The notion of size and maximal size is analogous to $\alpha$-packings and uniform $\alpha$-packings.  For instance, the following lemma is immediate from the definitions:

\begin{lemma}\label{packing}  If $\P$ is a uniform $\alpha$-packing of a tree $T$, and $a(P)$ obeys the weak Carleson condition
$$ |a(P)| \leq A |I_P| \hbox{ for all } P \in \P$$
for some $A > 0$, then $\| a\|_{\size^*(\P)} \leq A\alpha$.
\end{lemma}

If $\mu$ is a non-negative measure on the truncated upper half-space $(R^2_+)_M$, then it assigns a non-negative number to each Whitney box.  By the correspondence \eqref{correlation}, we can thus assign to each lacunary tile $P = P^+(I)$ a number $\mu(P)$ by the formula $\mu(P^+(I)) := \mu(Q^+(I))$.  We say that $\mu$ is a \emph{Carleson measure} if
$$ \| \mu \|_{\size^*(\P^+)} < \infty.$$
The reader may easily verify that this is equivalent to the usual formulation of a Carleson measure, namely that $\mu(Q(I)) \leq C |I|$ for some constant $C$.

\subsection{Wavelets, phase space projections, and BMO}

Let $P$ be a lacunary tile.  We define the \emph{(mother) Haar wavelet} $\phi_P$ to be the $L^2$-normalized function
$$ \phi_P := |I_P|^{-1/2} (\chi_{I_P^{l}} - \chi_{I_P^{r}})$$
where $I_P^{l}$ and $I_P^{right}$ are the left and right halves of $I_P$ respectively.  Similarly, if $P$ is a non-lacunary tile, we define the \emph{(father) Haar wavelet} $\phi_P$ by
$$ \phi_P := |I_P|^{-1/2} \chi_{I_P}.$$
Observe that these functions are normalized in $L^2$, and that $\phi_P$ and $\phi_{P'}$ are orthogonal whenever $P$ and $P'$ are disjoint\footnote{For a pair of lacunary tiles, this means that $I_P \neq I_{P'}$; for a pair of non-lacunary tiles, this means that $I_P$ and $I_{P'}$ are disjoint.  For lacunary $P$ and non-lacunary $P'$, this means that $I_{P'}$ is not a proper subset of $I_P$.}.

It is in fact possible to assign a function $\phi_P$ to every Heisenberg tile; these functions are known as \emph{Walsh wave packets}, see e.g. \cite{thiele-walsh} for a discussion.  These Walsh packets can then be used to efficiently decompose such operators as the (Walsh) bilinear Hilbert transform or the (Walsh) Carleson operator, just as the Haar wavelets can be used to decompose (dyadic) paraproducts or (dyadic) Calder\'on-Zygmund operators; see in particular the remarks after \eqref{tril}.  However, we shall not make any use of the Walsh wave packets for the results in this paper.

We define a \emph{(dyadic) test function} to be any finite linear combination of mother and father Haar wavelets $\phi_P$.  We use $\S$ to denote the space of all test functions, and $\S_0$ to denote the test functions with mean zero. For any dyadic interval $I$, we define $\S(I)$ to be the elements of $\S$ which are supported in $I$, and similarly define $\S_0(I)$.  Note that $\S(I)$ is only one dimension larger than $\S_0(I)$, and is in fact spanned by $\S_0(I)$ and any function $b_I \in \S(I)$ with non-zero mean.  This fact will be used much later on when we discuss local $T(b)$ theorems.

Since the mother Haar wavelets are orthonormal, we have the representation formula\footnote{The continuous version of this would be a Calder\'on reproducing formula such as $f(x) = \int 2t^2\Delta e^{t^2\Delta}f(x) \frac{dt}{t}$.}
$$ f = \sum_{P \in \P^+} \langle f, \phi_P \rangle \phi_P$$
for all $f \in \S_0$.

If $f \in \S$ and $T$ is any collection of disjoint tiles in $\P^+ \cup \P^0$ 
(i.e. the tiles in $T$ can be lacunary or non-lacunary), we define the 
phase space projection $\Pi_T$ by
$$ \Pi_T f := \sum_{P \in T} \langle f, \phi_P \rangle \phi_P.$$
This is an orthogonal projection from $L^2(\R)$ to the space of functions spanned by $\{ \phi_P: P \in T\}$.  For instance, we have
$\Pi_{P^0(I)} f = [f]_I \chi_I$
and
\be{piti}
\Pi_{\Tree(I)} f = (f - [f]_I) \chi_I.
\end{equation}
More generally, for any convex tree $T \subset \P^+$ and any $f \in \S$ we have
\be{diff}
\Pi_T f(x) = [f]_{J(x,T)} - [f]_{I(x,T)}
\end{equation}
for some intervals $I(x,T)$, $J(x,T)$ containing $x$; the exact choice of these intervals depends on $T$.  The formula \eqref{diff} can be derived by writing $T$ as a complete tree $\Tree(P_T)$ with some smaller complete trees removed, and then using \eqref{piti}.

If $f \in \S$, we define the \emph{wavelet transform} $Wf$ of $f$ to be the function
$$ Wf(P) :=  \langle f, \phi_P \rangle$$
defined on $\P^+$; this is an isometry\footnote{The continuous analogue of this in the upper half-plane with measure $\frac{dx dt}{t}$ would be the function $Q_t f(x)$, where $Q_t$ is a suitable cancellative averaging operator with wavelength $t$, e.g. $Q_t := t^2 \Delta e^{t^2 \Delta}$.} between $\S_0$ (endowed with the $L^2$ norm) and $l^2(\P^+)$.
In particular, the function $|Wf|^2$ maps $\P^+$ to $\R^+$, and so one can compute the size of $|Wf|^2$ on various collections of trees.  We observe in particular that
\be{w-size}
\| |Wf|^2 \|_{\size(T)} = \frac{1}{|I_T|} \| \Pi_T f \|_2^2 \leq \frac{1}{|I_T|} \int_{I_T} |f - [f]_{I_T}|^2
\leq \frac{1}{|I_T|} \int_{I_T} |f|^2
\end{equation}
for all trees $T$.

If $f \in \S$, we define the \emph{dyadic BMO norm} of $f$ by 
\be{bmo-def}
\| f \|_{BMO} := \| |Wf|^2 \|_{\size^*(\P^+)}^{1/2}
\end{equation}
The reader may easily verify that this definition corresponds to the usual ($L^2$-based) definition of dyadic BMO.  Note that the projections $\Pi_T$ 
defined earlier are bounded on $L^2$ and BMO (in fact they are contractions).

Thus the concepts of Carleson measure, BMO, and maximal size are essentially the same concept.  However, the concept of maximal size extends more easily to general families of tiles (not necessarily lacunary) than the other two notions.
(In particular, the notion of maximal size on a tree centered at an arbitrary
frequency $\xi_0$ is central to the boundedness of the bilinear Hilbert
transform, see \cite{laceyt1}).

\subsection{Mean}

We shall need a notion of \emph{mean} (or ``normalized mass''), which can be thought of as a ``non-lacunary'' variant of size.  Given any function $f$ on $\R$ and a tile $P \in \P^+$, we define
$$ \| f \|_{\mean(P)} := [|f|]_{I_P} = \frac{1}{|I_P|} \int_{I_P} |f|$$
and for any collection $\P \subset \P^+$ of lacunary tiles we define
$$ \| f \|_{\mean^*(\P)} := \sup_{P \in \P} \| f \|_{\mean(P)}.$$
Like the notion of $BMO$, the notion of mean has the scaling of $L^\infty$.  One can extend the notion of mean to arbitrary tiles; the function $f$ should then be replaced by a measure on phase space.  For instance, in applications to Carleson's theorem \cite{carleson} the notion of mean is applied to a measure of the form $\chi_E(x) \delta(\xi - N(x))$, where $N$ is an arbitrary function.  See \cite{fefferman}, \cite{laceyt-carleson}.

\section{The ``$L^\infty$ theory''}\label{jn-sec}

It is well known that the notion of maximal size or BMO can be thought of as a stable substitute for the $L^\infty$ norm, which is often ill-suited for applications.  In this section we develop the standard theory for this norm.  

We begin with a simple but very useful principle: to bound the maximal size of a collection of tiles, it suffices to do so outside of an $(1-\eta)$-packing of each tree in the collection. 

\begin{figure}[htbp]\label{tree-fig} \centering
\ \psfig{figure=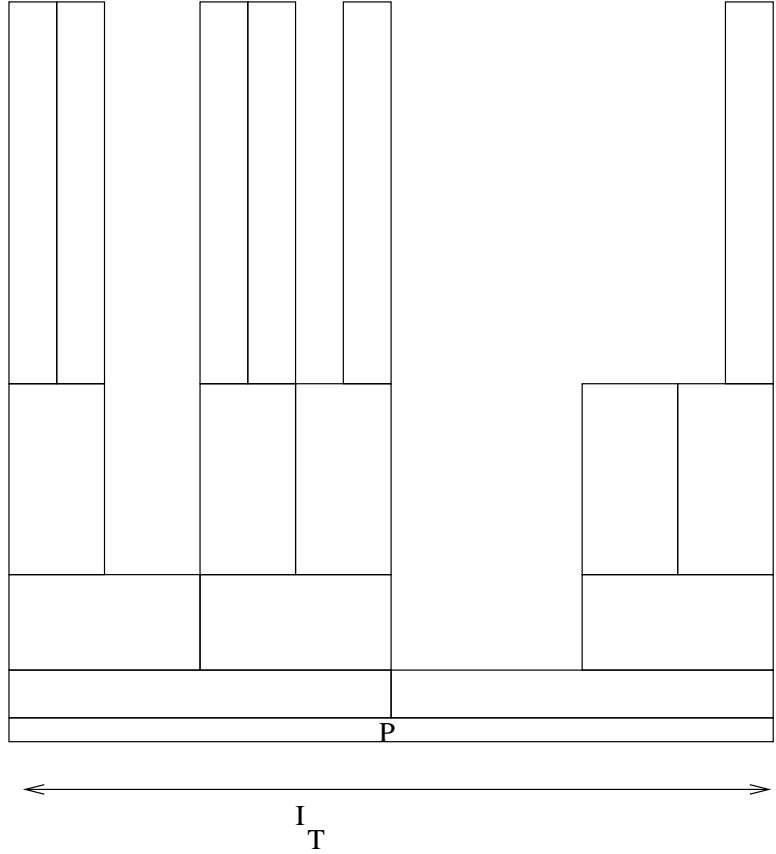}
\caption{A convex tree $T$ with top $P_T = P$.  Note how this tree can be thought of as a complete tree $\Tree(P)$ with some smaller complete trees removed.  In this particular tree, the tops of the trees removed form a $\frac{5}{8}$-packing of $\Tree(P)$, so that $\frac{3}{8}$ of the tree $T$ ``makes it all the way to the top''.  On the Carleson half-plane, this region resembles the portion of a Carleson box above a Lipschitz graph.}
\end{figure}

\begin{lemma}\label{tree-samp}  Suppose $\P$ is a collection of lacunary tiles, and let $a: \P \to \R^+$, $A > 0$, and $0 < \eta < 1$ be such that for every tree $T$ which is complete with respect to $\P$, one has 
$$ \| a \|_{\size(T \backslash \bigcup_{T' \in \T} T')} \leq A$$
for some collection $\T$ of trees in $T$ whose tops $\{ P_{T'}: T' \in \T \}$ are a $(1-\eta)$-packing of $T$.
Then we have
$$ \| a \|_{\size^*(\P)} \leq A/\eta.$$
\end{lemma}

\begin{proof}
Let $T$ be a tree in $\P$.  From hypothesis we have
\bas
\sum_{P \in T} a(P) &= \sum_{P \in T \backslash \bigcup_{T' \in \T} T'} a(P) + \sum_{T' \in \T} \sum_{P \in T'} a(P) \\
&\leq A |I_T| + \sum_{T' \in \T} \| a \|_{\size^*(\P)} |I_{T'}| \\
&\leq |I_T| (A + (1-\eta) \| a\|_{\size^*(\P)}).
\end{align*}
Dividing by $|I_T|$ and taking suprema of both sides we obtain
$$ \| a \|_{\size^*(\P)} \leq A + (1-\eta) \| a\|_{\size^*(\P)}$$
and the claim follows.
\end{proof}

An alternate (and perhaps more intuitive) proof of Lemma \ref{tree-samp} is to start with a tree $T$, estimate the ``good'' part $T \backslash \bigcup_{T' \in \T} T'$ of the tree, and then pass to the ``bad'' trees $T' \in \T$ and iterate this process until the tree is completely exhausted.  Since the geometric series $\sum_{n=0}^\infty (1-\eta)^n$ converges to $1/\eta$, the claim follows.

\begin{corollary}[Good-lambda characterization of maximal size]\label{good-lambda}  Let $\P$ be a collection of lacunary tiles, and let $a: \P \to \R^+$ be such that
$$ | \{ x \in I_T: \sum_{P \in T} a(P) \frac{\chi_{I_P}(x)}{|I_P|} \geq A \}| \leq (1 - \eta) |I_T|$$
for some $A > 0$, $0 < \eta < 1$ and all trees $T \subseteq \P$.  Then we have
$\| a \|_{\size^*(\P)} \leq \frac{A}{\eta}$.
\end{corollary}

This lemma goes back at least to Fritz John \cite{john}.  A partial converse can be obtained from Markov's inequality:
$$
| \{ x \in I_T: \sum_{P \in T} a(P) \frac{\chi_{I_P}(x)}{|I_P|} \geq A \} |
\leq \frac{1}{A}\| \sum_{P \in T} a(P) \frac{\chi_{I_P}}{|I_P|} \|_1
\leq \frac{1}{A} |I_T| \| a \|_{\size(T)}
\leq \frac{1}{A} |I_T| \| a \|_{\size^*(\P)}.$$

\begin{proof}
Let $T$ be any tree in $\P$.  Consider the set $\Q$ of all tiles $Q \in T$ such that
$$ \sum_{P \in T: Q \leq' P} a(P)/|I_P| \geq A$$
and such that $Q$ is maximal with respect to the ordering $<'$.
By assumption $\Q$ is a $(1-\eta)$-packing of $T$. 

By Lemma \ref{tree-samp} it suffices to show that
$\| a \|_{\size(T \backslash \bigcup_{Q \in \Q} \Tree(Q))} \leq A$,
or equivalently that
$$ \int_{I_T} \sum_{P \in T \backslash \bigcup_{Q \in \Q} \Tree(Q)} a(P) \frac{\chi_{I_P}(x)}{|I_P|}\ dx \leq A |I_T|.$$
But by construction of $\Q$, the integrand is bounded by $A$ for all $x \in I_T$, and the claim follows.
\end{proof}

A similar argument gives the well-known $L^p$ characterization of BMO:

\begin{corollary}\label{jn-char}
Let $0 < p < \infty$, and let $f \in \S$ be such that
$\frac{1}{|I|} \int_I |f - [f]_I|^p \lesssim 1$,
or equivalently that
$\int_I |\Pi_{\Tree(I)} f|^p \lesssim |I|$
for all dyadic intervals $I$.  Then
$$ \|f\|_{BMO} \lesssim 1.$$
(The implicit constants depend on $p$).
\end{corollary}

Applying this corollary with $p=1$ we obtain in particular that
$$ \| f \|_{BMO} \sim \sup_I \frac{1}{|I|} \int_I |f - [f]_I|;$$
in fact, this is sometimes taken as the definition of (dyadic) BMO.

\begin{proof}
We need to show that $|Wf|^2$ has bounded maximal size.  Let $T = \Tree(I)$ be any complete tree for some interval $I$.  By Lemma \ref{tree-samp}, it suffices to find a collection $\T$ of trees in $T$ whose tops are a $\frac{1}{2}$-packing of $T$ such that
\be{wf-bound}
\| |Wf|^2 \|_{\size(T \backslash \bigcup_{T' \in \T} T')} \lesssim 1
\end{equation}

First observe from \eqref{piti} that if $J \subset I$ and $x \in J$ then
$$ \Pi_{\Tree(I)} f(x) - [\Pi_{\Tree(I)} f]_J = \Pi_{\Tree(J)} f(x)$$
so from hypothesis we have
\be{pifj}
\int_J |\Pi_{\Tree(I)} f - [\Pi_{\Tree(I)} f]_J|^p \lesssim |J|.
\end{equation}
Now let $C_p$ be a large constant to be chosen later.  Let $\Q$ denote the tiles $Q \in \Tree(I)$ such that
$|[\Pi_{\Tree(I)} f]_Q| \geq C_p$
and that $Q$ is maximal with respect to $\leq'$.  
If $C_p$ is sufficiently large, we see from \eqref{pifj} 
$$ \int_{I_Q} |\Pi_{\Tree(I)} f|^p  \gtrsim C_p^p |I_Q|$$
and hence that $\Q$ is a $\frac{1}{4}$-packing of $T$.  In particular the collection $2\Q = \{ 2Q: Q \in \Q \}$ of parents of tiles in $\Q$ is a $\frac{1}{2}$-packing of $T$.

We set $\T := \{ \Tree(2Q): 2Q \in 2\Q \}$, and define
$$ F := \Pi_{T \backslash \bigcup_{T \in \T} T} f = \Pi_{\Tree(I)} f - \sum_{Q \in \Q} \Pi_{\Tree(2Q)} f.$$
Then we can rewrite the left-hand side of \eqref{wf-bound} as
$\frac{1}{|I|} \| F \|_2^2$.

By construction we see that $F$ is supported on $I$.  Since $[\Pi_{\Tree(2Q)} f]_{2Q} = 0$, we see that $F$ is constant on each $I_{2Q}$ and that
$$ \|F\|_{L^\infty(2Q)} = |[F]_{2Q}| = |[\Pi_{\Tree(I)} f]_{2Q}| \lesssim C_p.$$
If $x \in I$ is not in any of the $I_{2Q}$, then
$$ |F(x)| = |\Pi_{\Tree(I)} f(x)| \lesssim C_p$$
Thus we have $\|F\|_\infty \lesssim C_p$.  Combining this with the previous we obtain \eqref{wf-bound} as desired. 
\end{proof}

We now give the well-known converse to the above Lemma:

\begin{lemma}[John-Nirenberg inequality]\label{jn}  Let $I$ be a dyadic interval, and let $f \in \S_0(I)$ be real-valued.  Then we have
$$ \| f \|_p \lesssim (1+p) |I|^{1/p} \| f \|_{BMO}$$
for all $0 < p < \infty$ and
$$ | \{ x \in I: f(x) > 2n \|f\|_{BMO} \} | \leq 2^{-n+1} |I| \hbox{ for all } n \in \Z^+.$$
\end{lemma}

\begin{proof}
It suffices to prove the latter inequality, as the former easily follows.

We prove the claim by induction on $n$.  The claim is clear for $n =1$.  Now suppose that $n > 1$ and the claim has already been proven for $n-1$.

Fix $I$, $f$.  Let $\P$ denote those tiles $P$ in $\Tree(I)$ such that $[f]_P > 2 \|f\|_{BMO}$, and such that $\P$ is maximal with respect to $\leq'$.  For each $P$ we have 
$$ \int_{I_P} |f|^2 \geq |I_P| |[f]_P|^2 \geq 4 |I_P| \|f\|_{BMO}^2.$$
On the other hand, from \eqref{w-size}, \eqref{bmo-def} we have
$$ \int_I |f|^2 \leq |I| \| |Wf|^2 \|_{\size(\Tree(I))} \leq |I| \|f\|_{BMO}^2.$$
Thus $\P$ is a $\frac{1}{4}$-packing of $\Tree(I)$, so that the collection $2\P = \{ 2P: P \in \P\}$ of parents of tiles in $\P$ form a $\frac{1}{2}$-packing of $\Tree(I)$.

By construction we have $[f]_{2P} \leq 2 \|f\|_{BMO}$ for all $P \in \P$, and $f(x) \leq 2 \|f\|_{BMO}$ for all $x \not \in \bigcup_{P \in \P} I_{2P}$.  Thus 
$$ \{ x \in I: f(x) > 2n \|f\|_{BMO} \} \subseteq \bigcup_{2P \in 2\P}
\{ x \in I_{2P}: f-[f]_{2P}(x) > 2(n-1) \|f\|_{BMO} \}.$$
The claim then follows from the inductive hypothesis.
\end{proof}

\subsection{Chopping big trees into little trees, and extrapolation of Carleson measures}\label{chopping-sec}

Let $a: \P^+ \to \R^+$ be a function.  Suppose we have a convex lacunary tree $T_0$ with a large size, let's say
\be{size-t0}
\| a \|_{\size^*(T_0)} \leq C_0.
\end{equation}
Let $0 < \delta \leq C_0$ be a small number.  An obvious question to ask is whether one can decompose the large tree $T_0$ into small trees, each of which has size less than or equal to $\delta$.  This is clearly impossible, as the example of a singleton tree $T_0$ with large size demonstrates.  However, one can do the next best thing:

\begin{theorem}\label{tree-slice}  With the above assumptions, we have the disjoint partition
\be{t0-part}
T_0 = \bigcup_{T \in \T} T \cup \P
\end{equation}
where the trees $T$ in $\T$ are convex and satisfy
\be{tree-small}
\|a\|_{\size^*(T)} \leq \delta
\end{equation}
while the tiles $P \in \P$ obey the estimate
\be{weak-carleson}
a(P) \leq C_0 |I_P|.
\end{equation}
Furthermore, the tiles $\P$ and the tree tops $\{ P_T: T \in \T \}$ are both uniform $C(C_0,\delta)$-packings of $T_0$.
\end{theorem}

Note that Lemma \ref{packing} gives an easy converse to the above Theorem: if $T_0$ can be partitioned by \eqref{t0-part} with the above properties then $\|a\|_{\size^*(T_0)}$ is bounded (but by a much large constant than $C_0$).  Thus, if one is willing to ignore losses in constants, the above Theorem gives a complete characterization of trees of large size in terms of trees of small size.  As we shall see in this section, this theorem can be applied to give extrapolation lemma for Carleson measures, and seems likely to be useful in other contexts also.

A continuous parameter version of Theorem \ref{tree-slice} 
is at least implicit in
\cite{kato-extrap}, where, as here, it is used to prove the 
``Extrapolation Lemma for Carleson Measures" (see Corollary 
\ref{mu-perturb} below).
The latter, in its continuous parameter form, was then used to establish the 
``restricted version" of the Kato square root conjecture, for 
$L^{\infty}$ perturbations of real, symmetric, elliptic coefficient matrices.
The essential idea of the extrapolation method had previously been introduced
by J. Lewis in his work with M. Murray \cite{lm}
on the heat equation in non-cylindrical 
domains, and refined further by Lewis and one of the present authors
\cite{hl} in their work on parabolic and elliptic equations.  Similar ideas
had also appeared previously in the work of David and Semmes on uniform 
rectifiability:  indeed Theorem 4.5 is very closely related to the 
``Corona Decomposition" of \cite{ds}. 

Roughly speaking, in applications of the extrapolation method, the idea 
is first to show that
some ``scale-invariant estimate on cubes"  (like a Carleson measure estimate, 
a BMO estimate, or a reverse Holder or $A_{\infty}$ estimate for a weight) 
holds when some controlling Carleson measure is suitably small in a 
certain sense, 
which will be made precise in the sequel.  The term ``extrapolation" 
refers to the removal of the smallness restriction.  In that sense it is 
analogous to G. David's technique for bootstrapping the Lipschitz constant
(see, e.g., \cite{david}), although it is not clear whether there exists an 
explicit connection between the two methods.

In \cite{lm}, \cite{hl}, for example, the controlling Carleson measure was a 
condition on either the boundary of the domain, or on the coefficients of the
elliptic or parabolic operator, and one proved reverse Holder iequalities for 
the associated elliptic-harmonic or parabolic measures.  In particular,
in \cite{hl}, the authors give an alternative proof, via extrapolation,
of the main theorem of R. Fefferman, Kenig and Pipher
\cite{fkp}, in which the controlling Carleson measure is a condition on 
the disagreement between the coefficients of two elliptic (or parabolic)
operators, in the case that reverse Holder estimates are known to hold
for the 
elliptic-harmonic measures
associated to the first operator, and one wishes to 
prove such estimates for the second.

In \cite{kato-extrap},
the authors exploit the fact that proving Kato's square root estimate is 
equivalent (by ``$\CZ(1)$" type reasoning)
to proving that a certain positive measure in 
the upper half space is Carleson.  Here,
the controlling Carleson measure was the one associated
to the original, self-adjoint operator,
and the extrapolation technique was used to prove that the analogous measure, 
related to the square root estimate for the perturbed operator, was also 
Carleson.  It is in this setting that Corollary \ref{mu-perturb}, or rather its 
continuous parameter analogue, is directly applicable.

The proof we give here follows the approach in \cite{kato-extrap}.  At the end of this section we give an alternate proof, due to John Garnett, which gives better dependence on constants.

Before we give the rigorous proof, we first informally describe the idea of the argument.  Suppose the original tree $T_0$ has size $\| a \|_{\size(T_0)} = c$.  Then $0 \leq c \leq C_0$ by \eqref{size-t0}.  To create a tree of maximal size less than $\delta$, we start with $T_0$ and remove from it some sub-trees of size between $c + \delta/2$ and $c - \delta/2$, which we select by a straightforward stopping time argument.  It then remains to control the sub-trees that were removed.  By shrinking the trees slightly (putting the error into $\P$) we can assume that the trees have size either greater than $c+\delta/2$ or less than $c-\delta/2$ (so that the tree that remains must have size at most $\delta$).  We call the first type of tree ``heavy'' and the second type ``light''.  Because the original tree had size $c$, it cannot be the case that $I_{T_0}$ is covered by heavy sub-trees, and so a positive proportion of $I_{T_0}$ must be covered by light trees or by nothing.  We then pass to the light sub-trees and iterate this process, finding a positive proportion $I_{T_0}$ occupied by increasingly lighter sub-trees.  After about $O(C_0/\delta)$ steps, we must terminate, finding a positive proportion of $I_{T_0}$ which are not covered by any further sub-trees.  We then pass to the remaining portion of $I_{T_0}$ and all the heavy trees which have until now been neglected, and iterate once again; since we have replaced $I_{T_0}$ with a strictly smaller fraction of $I_{T_0}$, this procedure will converge geometrically to obtain the desired estimates.

We now prove Theorem \ref{tree-slice}.  We shall drop \eqref{weak-carleson} since it follows from \eqref{size-t0}.  In the spirit of Lemma \ref{tree-samp}, it will suffice to prove the apparently weaker 

\begin{theorem}\label{tree-slice-weak}  With the above assumptions, we can find a (possibly empty) collection $\T_{iterate}$ of disjoint convex trees in $T_0$ whose tops have disjoint spatial intervals and form a $(1-\eta)$-packing of $T_0$
for some $\eta = \eta(C_0,\delta) > 0$, such that we have the disjoint partition
\be{t0-part-weak}
T_0 = \bigcup_{T \in \T_{iterate}} T \cup \bigcup_{T \in \T} T \cup \P
\end{equation}
where the trees $T \in \T$ obey \eqref{tree-small}, and $\P$ and the tree tops of $\T$ are both uniform $C(C_0,\delta)$-packings of $T_0$.
\end{theorem}

Indeed, if Theorem \ref{tree-slice-weak} holds, then we can construct the collection in Theorem \ref{tree-slice} by starting with the partition \eqref{t0-part-weak}, and then taking each of the trees in $\T_{iterate}$ and breaking them up by a further application of Theorem \ref{tree-slice-weak}.  We continue on in this way until the original tree $\T_0$ is completely broken up into trees $T$ obeying \eqref{tree-small} and tiles $P$ obeying \eqref{weak-carleson}.  The fact that $\P$ and the tree tops of $\T$ are $C(C_0,\delta)$-packings of $T_0$ then follows from Theorem \ref{tree-slice-weak} and the fact that the geometric series $\sum_n (1-\eta)^n$ converges.  A similar argument can then be used to improve ``$C(C_0,\delta)$-packing'' to ``uniform $C(C_0,\delta)$-packing''. We omit the details.

\begin{proof} (of Theorem \ref{tree-slice-weak}).  Define the quantity $c$ by
\be{a-bound}
c := \| a \|_{\size(T_0)},
\end{equation}
thus $0 \leq c \leq C_0$.  We shall prove the theorem by induction on $c$.  Specifically, we fix $0 \leq c \leq C_0$ and assume that the theorem has already been proven in the case $\| a \|_{\size(T_0)} \leq c-\delta/2$.  Note that we only have to apply this induction a finite number of times (about $O(C_0/\delta)$) so we will be allowed to let the constants get worse with each induction step.

The main lemma used in the proof of the Theorem will be

\begin{lemma}\label{slice}  We can partition
\be{part}
T_0 = \bigcup_{T \in \T_{small}} T \cup \P_{buffer} \cup \bigcup_{T \in \T_{heavy}} T \cup \bigcup_{T \in \T_{light}} T
\end{equation}
where $\T_{small}$ is a collection of convex trees which all obey \eqref{tree-small} and whose tree tops are a uniform $4$-packing of $T_0$,
$\P_{buffer}$ is uniform $3$-packing of $T_0$, and $\T_{heavy}$, $\T_{light}$ are collections of disjoint convex sub-trees of $T_0$ which are complete with respect to $T_0$, and are such that we have the tree counting estimates
\begin{align}
\sum_{T \in \T_{heavy}} |I_T| + \sum_{T \in \T_{light}} |I_T| &\leq |I_{T_0}|\label{t-1}\\
\sum_{T \in \T_{heavy}} |I_T|  &\leq \frac{c}{c + \delta/2} |I_{T_0}|\label{t-2}
\end{align}
and the size bounds
\be{lightness}
\| a\|_{\size(T)} \leq c-\delta/2 \hbox{ for all } T \in \T_{light}.
\end{equation}
\end{lemma}

\begin{figure}[htbp]\label{extrap-fig} \centering
\ \psfig{figure=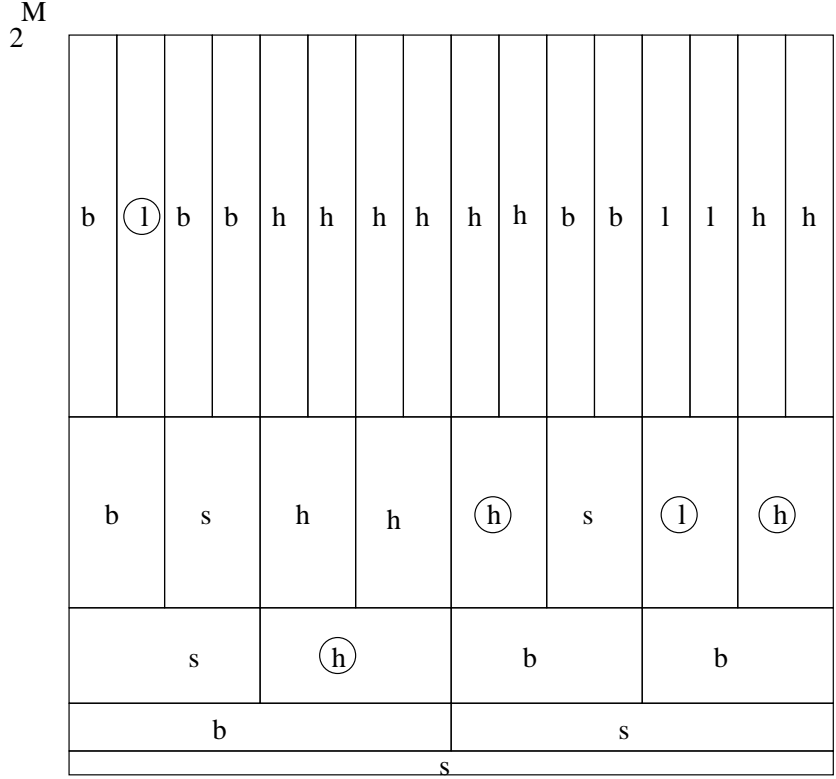}
\caption{A convex tree $T_0$ and its decomposition from Lemma \ref{slice}.  
The circled $h$ and $l$ tiles are the tops of maximal sub-trees of $T_0$ for 
which the size of $a$ fluctuates by at least $\delta/2$ from $c$; the
uncircled $h$ and $l$ tiles are the remaining tiles in those maximal
sub-trees.  $\T_{heavy}$ thus consists of the (three) $h$ trees while $\T_{light}$
consists of the (two) $l$ trees.  $\P_{buffer}$ consists of those remaining tiles
(labeled $b$) which lie just below a heavy or light tile, or are at the very top
of the phase plane.  The remaining tiles (labeled $s$) form the (three) small trees
$\T_{small}$.}
\end{figure}

\begin{proof}
Define $\T_{fluctuate}$ to be those sub-trees $T$ of $T_0$ 
which are complete with respect to $T_0$, such
that 
$| \| a \|_{\size(T)} - c | \geq \delta/2,$
and such that $T$ is maximal with respect to set inclusion and the above two properties.  Note that such trees are automatically convex.

By construction, none of the trees in $\T_{fluctuate}$ contain the top tile $P_{T_0}$.  We may subdivide\footnote{With reference to Figure \ref{extrap-fig}, $\T_{fluctuate}$ consists of the $h$ and $l$ trees, $T_1$ consists of the $s$ and $b$ tiles, and $T_2$ consists of just the $s$ tiles.} 
$$ \T_{fluctuate} = \T_{heavy} \cup \T_{light}$$
where $\T_{heavy}$ consists of those trees $T \in \T_{fluctuate}$ with
\be{heavy}
\| a \|_{\size(T)} \geq c + \delta/2
\end{equation}
and $\T_{light}$ consists of those trees $T \in \T_{fluctuate}$ with
$$ \| a \|_{\size(T)} \leq c - \delta/2.$$
The trees in $\T_{fluctuate}$ are disjoint, convex, and have disjoint spatial supports, so \eqref{t-1} holds.  On the other hand if one multiplies \eqref{heavy} by $|I_T|$ and sums over all $T \in \T_{heavy}$ one obtains
$$ |I_{T_0}| c = \sum_{P \in T_0} a(P)
\geq \sum_{P \in \bigcup_{T \in \T_{heavy}} T} a(P)
\geq (c + \delta/2) \sum_{T \in \T_{heavy}} |I_T|.$$
Dividing by $c + \delta/2$ we obtain \eqref{t-2}.

Let $T_1$ denote the convex tree
$T_1 := T_0 \backslash \bigcup_{T \in \T_{fluctuate}} T$
with top $P_{T_0}$.  Informally, $T_1$ represents the portion of $T_0$ below the fluctuating tiles.  The tree $T_1$ contains $P_{T_0}$ and is hence non-empty.  Let $\P_{buffer}$ denote the tiles\footnote{Here we are taking advantage of our decision to work in a finite model, where the tiles have a minimal width $2^{-M}$.  One can replicate this argument in the infinite setting but one has to treat the portion of $T_1$ which ``goes all the way to infinity'' separately.  See \cite{kato-extrap}.} 
$$ \P_{buffer} := \{ P \in T_1: P=2Q \hbox{ for some } Q \not \in T_1 \} \cup \{ P \in T_1: |I_P| = 2^{-M} \}.$$
In other words, $\P_{buffer}$ consists of those tiles in $T_1$ which touch the upper boundary of $T_1$ (which in particular may include the tiles of minimal width $|I_P| = 2^{-M}$).  Since the tiles $Q$ in the definition of $\P_{buffer}$ have disjoint spatial supports and $|I_P| = 2|I_Q|$ we see that
$\{ P \in T_1: P=2Q \hbox{ for some } Q \not \in T_1 \}$ is a uniform
2-packing of $T_1$.  Since $\{ P \in T_1: |I_P| = 2^{-M} \}$ is clearly
a uniform 1-packing of $T_1$, we thus see that $\P_{buffer}$ is a uniform
3-packing of $T_1$.

Let $T_2$ denote the (possibly empty) tree
$T_2 := T_1 \backslash \P_{buffer}$
with top $P_{T_0}$.  This tree is not necessarily convex, however we shall invoke the following lemma to split it into convex trees.

\begin{lemma}\label{convexify}  Let $T$ be a convex tree, and let $\P \subset T$ be a uniform $\alpha$-packing of $T$ for some $\alpha > 0$.  Then $T \backslash \P$ can be partitioned into
$$ T \backslash \P = \bigcup_{T' \in \T} T'$$
where $\T$ is a collection of convex trees $T'$ whose tops $\{ P_T: T \in \T \}$ form a uniform $(\alpha+1)$-packing of $T$. 
\end{lemma}

\begin{proof}
Let $\Q$ denote those dyadic intervals $Q \subset I_T$ such that $Q \in T \backslash \P$ and $2Q \not \in T \backslash \P$.  For any $Q \in \Q$, we see that either $Q = P_T$ or $2Q \in \P$.  Since $\P$ is a uniform $\alpha$-packing, this implies that $\Q$ is a uniform $(\alpha+1)$-packing. 

For each $Q \in \Q$, define the convex tree $T_Q$ with top $Q$ by
$$ T_Q := \{ P \in T \backslash \P: P \leq' Q, \hbox{ and there does not exist } Q' \in \Q \hbox{ such that } P <' Q' <' Q \}.$$
If we then set
$\T := \{ T_Q: Q \in \Q \}$ we see that the Lemma follows.
\end{proof}

By Lemma \ref{convexify} we may write $T_2 = \bigcup_{T \in \T_{small}} T$ where the trees in $\T_{small}$ are distinct and the tree tops of $\T_{small}$ are a uniform 4-packing of $\T$.

We now verify that each tree $T \in \T_{small}$ obeys \eqref{tree-small}.
It suffices to show that
\be{subtract}
\sum_{P \in \Tree(I) \cap T} a(P) \leq \delta |I|
\end{equation}
for all $I \subseteq I_{T_0}$.  

Fix $I$.  The idea is to write $\Tree(I) \cap T$ as the difference of trees,
each of which has size $c + O(\delta)$.

We may assume that $P^+(I) \in T$ since the claim is trivial otherwise.  We observe that
$$ \Tree(I) \cap T = (\Tree(I) \cap T_0) \backslash \bigcup_{J \in \J} (\Tree(J) \cap T_0)$$
where $\J$ consists of those intervals $J \subseteq I$ such that $P^+(J) \not \in T$, and which are maximal with respect to this property.

The tile $P^+(I)$ is in $T$ and hence in $T_1$.  By construction of $T_1$, we thus have 
$$\sum_{P \in \Tree(I) \cap T_0} a(P) = |I| \| a\|_{\size(\Tree(I) \cap T_0)} 
\leq |I| (c + \delta/2 )$$
(since otherwise $\Tree(I) \cap T_0$ would belong to $\T_{heavy}$, a contradiction).
Similarly, for every $J \in \J$, the tile $\P^+(J)$ is contained in $T_1$ (otherwise $P^+(2J)$ would be both in $T$ and in $\P_{buffer}$, a contradiction), so 
$$\sum_{P \in \Tree(J) \cap T_0} a(P) 
\geq |J| (c - \delta/2 ).$$
By the construction of $\J$, the intervals $J$ in $\J$ partition $I$, thus
$$\sum_{J \in \J} \sum_{P \in \Tree(J) \cap T_0} a(P) 
\geq |I| (c - \delta/2 ).$$
Subtracting this from the previous we obtain \eqref{subtract} as desired.
\end{proof}

We apply the above lemma and place $\P_{buffer}$ into $\P$, $\T_{small}$ into $\T$, and $\T_{heavy}$ into $\T_{iterate}$.  For the remaining trees $\T_{light}$ we use the induction hypothesis, which splits each of the trees in $\T_{light}$ into $\T_{iterate}$, $\T$, and $\P$.  All the desired conclusions of Theorem \ref{tree-slice-weak} are easily verified except perhaps for the claim that the tops of $\T_{iterate}$ form a $(1-\eta)$-packing of $T$, or in other words
$$ \sum_{T \in \T_{iterate}} |I_T| \leq (1-\eta) |I_{T_0}|.$$
To prove this inequality, note that the trees in $\T_{heavy}$ contribute
$\sum_{T \in \T_{heavy}} |I_T|$
to the left-hand side, while from the induction hypothesis the trees in $\T_{light}$ contribute at most
$(1 - \eta) \sum_{T \in \T_{light}} |I_T|$
for some $\eta > 0$.
There are no other contributions. The claim then follows from \eqref{t-1}, \eqref{t-2} (reducing the value of $\eta$ as necessary).
\end{proof}

The constants $C(C_0,\delta)$ given by this argument are about $(C_0/\delta)^{C C_0/\delta}$.  This bound can be improved substantially; see below.
 
The following corollary allows one to use one Carleson measure $\mu$ to prove the Carleson measure property of a related measure $\mu'$.  It is the dyadic version of 
an extrapolation lemma in \cite{kato-extrap}, which in turn is based on ideas in \cite{lm}, \cite{hl}.

\begin{corollary}[Extrapolation of Carleson measures]\label{mu-perturb}  Let $\mu: \P^+ \to \R^+$ have bounded maximal size and let $\delta > 0$.  Let $\mu'$ be a non-negative measure on $\R^2_+$ obeying the ``weak Carleson condition''
$$ \mu'(P) \leq C_1 |I_P| \hbox{ for all } P \in \P^+$$
and such that $\| \mu' \|_{\size(T)} \leq C_2$
for all convex trees $T$ such that $\| \mu \|_{\size^*(T)} \leq \delta$.
Then $\mu'$ also has bounded maximal size:
$$ \| \mu' \|_{\size^*(\P^+)} \leq C(\| \mu\|_{\size^*(\P^+)}, \delta) (C_1 + C_2).$$
\end{corollary}

\begin{proof}
Let $T_0$ be any convex tree.  We need to show that
$$ \| \mu' \|_{\size(T_0)} \leq C(\| \mu\|_{\size^*(\P^+)}, \delta) (C_1 + C_2).$$
By Theorem \ref{tree-slice}, we can partition $T_0 = \bigcup_{T \in \T} T \cup \P$
where
$\| \mu \|_{\size^*(T)} \leq \delta$ for all $T \in \T$,
and
$$ \sum_{T \in \T} |I_T| + \sum_{P \in \P} |I_P| \leq C(\| \mu\|_{\size^*(\P^+)}, \delta)  |I_{T_0}|.$$
From this and assumptions on $\mu'$ we see that
$$ \sum_{P \in T_0} \mu'(P) = \sum_{T \in \T} \sum_{P \in T} \mu'(P) + \sum_{P \in \P}\mu'(P) \leq C(\| \mu\|_{\size^*(\P^+)}, \delta) C_2 |I_{T_0}|
+ C(\| \mu\|_{\size^*(\P^+)}, \delta) C_1 |I_{T_0}|$$
and the claim follows.
\end{proof}

As mentioned earlier, this lemma has applications to the Kato problem.  In \cite{kato-extrap}, this lemma was used to establish
a restricted version of Kato's conjecture, for perturbations of real, symmetric
coefficient matrices.  In that case, $\mu$ was a Carleson measure which
controlled the original operator, and $\mu'$ was the analogous measure
controlling the perturbed operator.  The point was to establish that $\mu'$ was
also a Carleson measure.  We remark that the fact that the final bound on $\mu'$ was linear in $C_2$ was crucial to this application.

It is possible to eliminate the weak Carleson condition by allowing the tree measured by $\mu'$ to be a little larger than the tree measured by $\mu$, but we will not pursue this type of generalization here.

\subsection{An alternate argument}

In this section we give an alternate proof of Theorem \ref{tree-slice}, due to John Garnett (personal communication).  The idea of this argument is similar
to some arguments in \cite{bcgj}.

Fix $T_0$, $a$.   We first observe that it suffices to prove the theorem under the additional ``weak Carleson'' assumption
\be{a-small}
a(P) \leq \frac{\delta}{2} |I_P| \hbox{ for all } P \in T_0.
\end{equation}
To see this, suppose that we are in the general case when \eqref{a-small} need not hold.  We set $\P$ to be the set of tiles where \eqref{a-small} fails:
\be{p-big}
\P := \{ P \in T_0: a(P) \geq \frac{\delta}{2} |I_P| \}.
\end{equation}
From \eqref{size-t0} we see that $\P$ is a uniform $2C_0/\delta$-packing of $T_0$ (cf. Lemma \ref{packing}).  By Lemma \ref{convexify} we thus see that we can split $T_0 \backslash \P$ into a collection of disjoint convex subtrees of $T_0$, whose tops form a uniform $2C_0/\delta+1$-packing of $T_0$.  On each such subtree \eqref{a-small} holds.  Thus if we apply Theorem \ref{tree-slice} to each sub-tree and then combine all the decompositions, we obtain the desired decomposition for the original tree $T_0$ (with the constants $C(C_0,\delta)$ worsened by a factor of $2C_0/\delta+1$).

Henceforth we assume \eqref{a-small}.  Under this assumption we will not need
$\P$ any more, and will set it equal to the empty set.  

We can assume without loss of generality that $T_0$ is a complete tree, since if $T_0$ is incomplete then one can replace $T_0$ by its completion, and extend $a$ by zero; note that the intersection of two convex trees is convex, so one does not lose convexity when one restricts back to $T_0$.

We can now make the technical assumption that the minimal tiles have large coefficient:
\be{end-a}
a(P) = \frac{\delta}{2} |I_P| \hbox{ whenever } P \in T_0 \hbox{ and } |I_P| = 2^{-M}.
\end{equation}
This is because in the general case one can simply increase $a(P)$ for these tiles to equal $\frac{\delta}{2} |I_P|$; observe that this only increases $\|a\|_{\size^*(T)}$ by at most $\delta/2$, so the claim follows by redefining $C_0$ as necessary.  This technical assumption is needed to make sure that a certain stopping argument always halts before it reaches the smallest scale.

We now use the greedy algorithm to select a subtree $T$ of $T_0$ of size roughly comparable to $\delta$:

\begin{lemma}\label{td}  Let $T_0$ be a complete tree, and let $a: T_0 \to \R^+$ obey \eqref{a-small} and \eqref{end-a}.  Then there exists a convex subtree $T \subseteq T_0$ of $T_0$ with $P_T = P_{T_0}$ such that
$$ \frac{\delta}{2} \leq \| a \|_{\size(T)} \leq \| a \|_{\size^*(T)} \leq \delta$$
\end{lemma}

\begin{proof}
Consider the top tile $P_{T_0}$.  If $a(P_{T_0}) = \frac{\delta}{2} |I_{T_0}|$ then we can set $T$ to be the singleton tree $\{ P_{T_0} \}$, so we can
assume that $a(P_{T_0})$ is strictly less than $\frac{\delta}{2} |I_{T_0}|$.

Let $\T$ be the class of all convex subtrees $T' \subseteq T_0$ with top $P_{T'} = P_{T_0}$ such that 
\be{ads}
| a \|_{\size^*(T')} < \delta/2.
\end{equation}
  By the previous paragraph, this class $\T$ is non-empty; also, by \eqref{end-a}, this class cannot contain any tiles with the minimal width $2^{-M}$.

Let $T_*$ be a tree in $\T$ which is maximal with respect to set inclusion.  Let $\P$ denote the set of tiles $P$ in $T_0 \backslash T_*$ such that $2P \in T_*$; these are the tiles which lie just above $T_*$.  Since $T_0$ is complete and $T_*$ does not contain any tiles of minimal width, we see that the spatial intervals $\{ I_P: P \in \P \}$ partition $I_{T_0}$.  In particular the $\P$ are a uniform 1-packing of $T_0$.  

Set $T := T_* \cup \P$; this is clearly a convex sub-tree of $T_0$ with top $P_T = P_{T_0}$. From \eqref{ads}, \eqref{a-small}, and the uniform 1-packing property of $\P$ we see that $\| a \|_{\size^*(T)} \leq \delta$ (cf. Lemma \ref{packing}).  It thus remains to show the lower bound on size, i.e.
\be{jg}
 \sum_{P \in T} a(P) \geq \frac{\delta}{2} |I_{T_0}|.
\end{equation}

Call a tile $Q \in T$ \emph{heavy} if $\| a \|_{\size(\Tree(Q) \cap T)} \geq \frac{\delta}{2}$, or in other words
\be{heavy-def}
\sum_{P \in T: P \leq' Q} a(P) \geq \frac{\delta}{2} |I_Q|.
\end{equation}
Observe that for every tile $P \in \P$, there must exist a heavy tile $Q \in T$ such that $P \leq' Q$, since otherwise one could add $P$ to the tree $T$ while retaining the property \eqref{ads}, contradicting the maximality of $T$.  

Let $\Q$ denote the set of heavy tiles $Q$ in $T$ which are maximal with respect to the ordering $\leq'$.  By the previous paragraph we see that the spatial intervals of $\Q$ partition $T_0$.  
If one adds up \eqref{heavy-def} for all such tiles one obtains \eqref{jg}.  The proof of the lemma is now complete.
\end{proof}

When one removes $T$ from the complete tree $T_0$ we obtain a union of disjoint complete tree, which are of course smaller than $T_0$ but still obey \eqref{a-small} and \eqref{end-a}.  Thus we can iterate the above lemma to obtain

\begin{corollary}
Let $T_0$ be a complete tree, and let $a: T_0 \to \R^+$ obey \eqref{a-small} and \eqref{end-a}.  Then we can partition $T_0 = \bigcup_{T \in \T} T$, where $\T$ is a collection of disjoint convex trees $T$ such that
$$ \frac{\delta}{2} \leq \| a \|_{\size(T)} \leq \| a \|_{\size^*(T)} \leq \delta$$
for all $T \in \T$.
\end{corollary}

We apply the above Corollary to the tree $T_0$ in the Theorem, obtaining the collection $\T$ of trees.  We now claim that the tops $\{P_T: T \in \T\}$ of these trees form a uniform $2C_0/\delta$-packing of $T_0$.  Indeed, for any dyadic interval $J \subseteq I_{T_0}$ we see that
\bas
\sum_{T \in \T: I_T \subseteq J}  \frac{\delta}{2} |I_T|
&\leq \sum_{T \in \T: I_T \subseteq J} \sum_{P \in T} a(P) \\
&\leq \sum_{P \in T_0: I_P \subseteq J} a(P) \\
&\leq \| a \|_{\size^*(T_0)} |J|\\
&\leq C_0 |J|.
\end{align*}
This completes the proof of Theorem \ref{tree-slice} (with $\P$ empty).
\endprf

Observe that the constants obtained in this manner are significantly superior to the previous argument, being polynomial in $C_0/\delta$ instead of exponential.

\section{The ``$L^p$ theory''}\label{tree-sec}

In the previous section we established some estimates in the case when the maximal size was bounded; this can be thought of as the ``$L^\infty$ theory'' of maximal size.  Now we study what happens when our collection of tiles does not have a good maximal size bound.  In this case we can subdivide the collection into disjoint trees, such that the size of each of the trees is under control:

\begin{lemma}[``Calder\'on-Zygmund decomposition for size'']\label{tree-select}
Let $n \in Z$, $\P_n$ be a convex collection of lacunary tiles, and let $a: \P_n \to \R^+$ be a function such that
$\| a \|_{\size^*(\P_n)} \leq 2^n$.
Then there exists a disjoint partition
\be{disj}
\P_n = \bigcup_{T \in \T_n} T \cup \P_{n-1}
\end{equation}
where $\P_{n-1}$ is a convex collection of tiles such that
\be{small-size}
\| a \|_{\size^*(\P_{n-1})} \leq 2^{n-1}
\end{equation}
and $\T_n$ is a collection of convex trees $T$ with disjoint spatial intervals $I_T$ such that
\be{med-size}
\| a \|_{\size(T)} \sim \| a \|_{\size^*(T)} \sim 2^n
\end{equation}
for all $T \in \T$.  

In the particular case that $a = |Wf|^2$ for some $f \in \S_0$, we then have
\be{pit-bmo}
\| \Pi_T f \|_{BMO} \sim 2^{n/2}
\end{equation}
and 
\be{pit-p}
\| \Pi_T f \|_p \sim |I_T|^{1/p} 2^{n/2}
\end{equation}
for all $0 < p < \infty$ (with the implicit constant depending on $p$).
\end{lemma}

\begin{proof}
We set $\T_n$ to be the collection of all trees $T \subset \P_n$ such that
$\| a \|_{\size(T)} \geq 2^{n-1}$,
and are maximal with respect to set inclusion.  Clearly these trees are disjoint (otherwise the union of the two trees would also qualify, and contradict maximality).  They are also complete with respect to $\P_n$ (i.e. $T = \Tree(I_T) \cap \P_n$) and thus convex.  One then sets
\be{pn1-def}
\P_{n-1} := \P_n \backslash \bigcup_{T \in \T_n} T.
\end{equation}
The properties \eqref{disj}, \eqref{small-size}, \eqref{med-size} are easily verified.  To prove the last two properties, observe from \eqref{med-size} that
$$ \| \Pi_T f \|_{BMO} = \| |Wf|^2 \|_{\size^*(T)}^{1/2} \sim 2^{n/2}$$
and
$$ \| \Pi_T f \|_2 = |I_T|^{1/2} \| |Wf|^2 \|_{\size(T)}^{1/2} \sim 2^{n/2} |I_T|^{1/2}.$$
The claim then follows from the John-Nirenberg inequality (Lemma \ref{jn}) and H\"older's inequality.
\end{proof}

The above Lemma should be compared with the standard Calder\'on-Zygmund decomposition, which if given a function $f$ with $\|f\|_\infty \leq 2^n$, will subdivide $f = g + \sum_I b_I$ where $\|g\|_\infty \leq 2^{n-1}$, the $b_I$ are supported on disjoint intervals $I$ and have mean zero and $\| b_I \|_p \sim 2^n |I|^{1/p}$ for all $0 < p < \infty$.  The trees in $\T_n$ are the analogues of the intervals $I$, and can be thought of as the region of phase space where $a$ (or $f$) ``has size $\sim 2^n$''.  Lemma \ref{tree-select} can also be thought of as a sort of BMO version of the Chebyshev's inequality
\be{cheby}
|\{ x: |f(x)| \gtrsim 2^{n/2} \}| \lesssim 2^{-np/2} \|f\|_p^p.
\end{equation}
Indeed, if $f$, $n$, $\T_n$ is as in the above lemma; then by \eqref{pit-p} and the disjointness of the $I_T$ we have
$$ \| f\|_p^p \gtrsim \sum_{T \in \T_n} \| \Pi_T f \|_p^p \sim \sum_{T \in \T_n} |I_T| 2^{np/2}$$
and hence
\be{cheb0}
|\bigcup_{T \in \T_n} I_T| \lesssim 2^{-np/2} \| f\|_p^p
\end{equation}
(compare with \eqref{cheby}, \eqref{pit-bmo}, \eqref{small-size}). 

In practice one iterates the above lemma, starting with a large $n$ and
decrementing $n$ repeatedly, thus decomposing $\P^+$ into trees of various
sizes (plus a remainder of size $0$).

We have a similar selection lemma for mean, which can be thought of
as the analogue of the previous lemma for the non-lacunary tiles $\P^0$.

\begin{lemma}[``Calder\'on-Zygmund decomposition for mean'']\label{mean-select}
Let $n \in Z$, $\P_n$ be a convex collection of lacunary tiles, and let $f \in \S$ such that
$\| f \|_{\mean^*(\P_n)} \leq 2^n$.
Then there exists a disjoint partition
\be{disj-mean}
\P_n = \bigcup_{T \in \T_n} T \cup \P_{n-1}
\end{equation}
where $\P_{n-1}$ is a convex collection of tiles such that
\be{small-mean}
\| f \|_{\mean^*(\P_{n-1})} \leq 2^{n-1}
\end{equation}
and $\T_n$ is a collection of convex trees $T$ with disjoint spatial intervals $I_T$ such that
\be{med-mean}
\| f \|_{\mean(P_T)} \sim \| f \|_{\mean^*(T)} \sim 2^n
\end{equation}
for all $T \in \T$.  In particular, we have
\be{cheb}
\sum_{T \in \T_n} |I_T| \lesssim 2^{-n} \int_{|f| \gtrsim 2^n} |f|
\lesssim 2^{-np} \|f\|_p^p
\end{equation}
for any $1 \leq p < \infty$ (with the implicit constant depending on $p$).
\end{lemma}

\begin{proof}
We define $\T_n$ to be the set of all trees $T$ such that $\|f\|_{\mean^*(T)} \geq 2^{n-1}$, and are maximal with respect to set inclusion, and then define $\P_{n-1}$ by \eqref{pn1-def}.  As before the trees $T$ are convex and complete with respect to $\P_n$.  The properties \eqref{disj-mean}, \eqref{small-mean}, \eqref{med-mean} are easily verified.  By \eqref{med-mean} we have
$\int_{x \in I_T: |f(x)| \gtrsim 2^n} |f| \gtrsim 2^n |I_T|$,
and \eqref{cheb} follows from the disjointness of the $I_T$.
\end{proof}

In later sections we give some applications of the above machinery.  These applications will all have a similar flavor, in that they follow the following broad strategy:

\begin{itemize}
\item Begin with a sum over a collection of lacunary tiles.
\item Use Lemma \ref{tree-select} and/or Lemma \ref{mean-select} to extract disjoint trees in this collection of a certain size (plus a remainder of size 0, which is usually trivial to handle).
\item Estimate the contribution of each tree in terms of the width $|I_T|$ of the tree and the size and/or mean of the tree.
\item Estimate the total width $\sum_T |I_T|$ of the trees (using such estimates as \eqref{pit-bmo}, \eqref{pit-p}, \eqref{cheb0}, \eqref{cheb}).
\item Remove these trees from the collection, and repeat the above steps until the collection has been exhausted.
\item Sum up.
\end{itemize}

This type of argument is fundamental to the general theory of tiles, as can be seen in the work on the bilinear Hilbert transform and Carleson's theorem (see e.g. \cite{fefferman}, \cite{laceyt1}, \cite{laceyt-carleson}, \cite{mtt:biest}).  Apart from several technical details, the main differences between the arguments in those papers and the ones here are to eliminate the word ``lacunary'' from the above strategy, and replace the above Lemmata by more sophisticated tree selection algorithms.

\subsection{Example: Atomic decomposition of dyadic $H^p$}\label{dyad-sec}

Let $0 < p \leq 1$; implicit constants will be allowed to depend on $p$.    In this section we reprove the standard atomic decomposition of $H^p$, but we first need some notation.

For any $f \in \S$, we define the \emph{dyadic Littlewood-Paley  square function} $Sf$ to be the vector-valued function
$$ Sf(x) := ( \Pi_P f(x) )_{P \in \P^+}
= ( \langle f, \phi_P \rangle \phi_P(x) )_{P \in \P^+}
= ( Wf(P) \phi_P(x) )_{P \in \P^+}$$
taking values in $l^2(\P^+)$.  Note that
\be{s-mag}
|Sf(x)| = (\sum_{P \in \P^+} |Wf(P)|^2 \frac{\chi_{I_P}(x)}{|I_P|})^{1/2}.
\end{equation}
The adjoint operator $S^*$ takes $l^2(\P^+)$-valued functions to scalar-valued
functions, and is given by the formula
$$
S^* (f_P)_{P \in \P^+} = \sum_{P \in \P^+} \Pi_P f_P.$$
In particular $S^* S$ is the identity on $\S_0$.  Also observe\footnote{Here we are allowing $W$ to act on vector-valued functions in the obvious manner, i.e. if $f = (f_1, \ldots, f_n)$ is a vector, then $Wf(P) := (Wf_1(P), \ldots, Wf_n(P))$.  Similarly we can define vector-valued BMO, etc.} that $|WSf(P)| = |Wf(P)|$ for all $P \in \P^+$, so in particular we have 
\be{isom}
\| Sf\|_2 = \|f\|_2; \quad \| Sf \|_{BMO} = \|f\|_{BMO}
\end{equation}
on $\S_0$ and $\S$ respectively. 

We also define the \emph{cancellative dyadic maximal function} $\tilde Mf$ by
$$ \tilde Mf(x) := \sup_{P \in \P^0} |\Pi_P f(x)| = \sup_{I: x \in I} \frac{1}{|I|} |\int_I f| = \sup_{I: x \in I} |[f]_I|;$$
this operator should not be confused with the dyadic Hardy-Littlewood maximal  function $Mf := \tilde M|f|$.  From \eqref{diff} we have
\be{pit-max}
|\Pi_T f(x)| \leq 2 \tilde Mf(x)
\end{equation}
whenever $T$ is a convex tree in $\P^+$.  

If $\Tree(I)$ is a complete tree, we define a \emph{(dyadic) $H^p$ atom on $\Tree(I)$} to be a function $a \in \S_0(I)$ such that $\| a \|_2 \leq |I|^{1/2 - 1/p}$.
Equivalently, $a \in \S_0$ is an $H^p$ atom on $\Tree(I)$ if and only if the wavelet transform $Wa$ of $a$ is supported on $\Tree(I)$ and $\| |Wa|^2 \|_{\size(\Tree(I))} \leq |I|^{-2/p}$; this is because of \eqref{w-size}.

In this section we show

\begin{theorem}[Equivalent definitions of $H^p$]\label{equiv}  Let $f \in \S_0$ and $0 < p \leq 1$.  Then the following statements are equivalent.
\begin{enumerate}
\renewcommand{\labelenumi}{(\roman{enumi})}
\item $\| Sf \|_p \lesssim 1.$
\item $\| \tilde Mf \|_p \lesssim 1.$
\item There exists a collection $\I$ of dyadic intervals, and to each $I \in \I$ there exists a non-negative number $c_I$ and an $H^p$ atom $a_I$ on $\Tree(I)$ such that $f = \sum_I c_I a_I$ and $\sum_I c_I^p \lesssim 1$.
\end{enumerate}
\end{theorem}

\begin{proof}
We first show that (iii) implies (i) and (ii).  From the quasi-triangle inequality
$$ \|f+g\|_p^p \leq \|f\|_p^p + \|g\|_p^p$$
we see that it suffices to verify this on atoms, i.e. to show that
$\| Sa_I \|_p, \|\tilde Ma_I \|_p \lesssim 1$
whenever $a_I$ is a $H^p$ atom on $\Tree(I)$.

Fix $I$, $a$.  By construction $\tilde M a_I$ and $S a_I$ are supported on $I$, so by H\"older it suffices to show
$\| Sa_I \|_2, \|\tilde Ma_I \|_2 \lesssim |I|^{1/2 - 1/p}$.
But this follows from the $L^2$ normalization of $a_I$ and the fact that $S$, $\tilde M$ are bounded on $L^2$.

It remains to show that either one of (i) or (ii) are enough to imply (iii).  Let $f$ be any element of $\S_0$, thus $f = \sum_{P \in \P^+} Wf(P) \phi_P$.

Set $a := |Wf|^2$.  We apply Lemma \ref{tree-select} repeatedly, starting with a sufficiently large $n$ and setting $\P_n := \P^+$, and then decrementing $n$ indefinitely.  Eventually one obtains a partition
$$ \P^+ = \bigcup_{n \in \Z} \bigcup_{T \in \T_n} T \cup \P_{-\infty}$$
where the $\T_n$ are as in Lemma \ref{tree-select}, and 
$\| |Wf|^2 \|_{\size^*(\P_{-\infty})} = 0$.
Thus $Wf$ vanishes on $\P_{-\infty}$, and only a finite number of $\T_n$ are non-empty.  We thus have
$f = \sum_{n \in \Z} \sum_{T \in \T_n} \Pi_T f$.
If we then set
$c_{I_T} := 2^{n/2} |I_T|^{1/p}$ and $a_{I_T} := \Pi_T f / c_{I_T}$
then we have
$$ f = \sum_{n \in \Z} \sum_{T \in \T_n} c_{I_T} a_{I_T}.$$

By \eqref{pit-p} (with $p=2$) we see that each $a_{I_T}$ is an $H^p$ atom.  To show (iii) it thus remains to show that
$$
\sum_n \sum_{T \in \T_n} c_{I_T}^p = \sum_n 2^{np/2} \sum_{T \in \T_n} |I_T| \lesssim 1.
$$

First suppose that (i) holds.  For each $n$ and each $T \in \T_n$, we see that
$$ \int_{I_T} |S \Pi_T f|^q
\lesssim |I_T| \| S \Pi_T f\|_{BMO}^q
= |I_T| \| \Pi_T f \|_{BMO}^q \sim 2^{nq/2} |I_T|$$
for all $2 < q < \infty$ by the John-Nirenberg inequality (Lemma \ref{jn}), \eqref{isom}, and \eqref{pit-bmo}.  
Also, we have
$$ \int_{I_T} |S \Pi_T f|^2 = \int_{I_T} |\Pi_T f|^2
\sim 2^{n} |I_T|.$$
By H\"older we thus have
$$
\int_{I_T} |Sf|^r \geq \int_{I_T} |S \Pi_T f|^r \gtrsim 2^{nr/2} |I_T|.$$
for any $0 < r < p$, with the implicit constant depending on $r$.  This clearly implies
$$ \int_{x \in I_T: |Sf(x)| \gtrsim 2^{n/2}} |Sf|^r \gtrsim 2^{nr/2} |I_T|.$$
Summing over all $T \in \T_n$ and using the disjointness of the $I_T$ we obtain
$$ \int_{ |Sf(x)| \gtrsim 2^{n/2}} |Sf|^r \gtrsim 2^{nr/2} \sum_{T \in I_T} |I_T|.$$
Multiplying by $2^{n(p-r)/2}$ and summing over $n$ we obtain
$$ \int \sum_{n: |Sf(x)| \gtrsim 2^{n/2}} 2^{n(p-r)/2} |Sf|^r \gtrsim \sum_n 2^{np/2} \sum_{T \in I_T} |I_T|.$$
Since the left-hand side is comparable to $\| Sf \|_p^p$, the claim follows.

Now suppose instead that (ii) holds.  By \eqref{pit-p}, \eqref{pit-max} we have
$
\int_{I_T} |\tilde Mf|^r \gtrsim 2^{nr/2} |I_T|$
for all $0 < r < p$.  Now we argue as with $Sf$ to obtain (iii) from (ii).
\end{proof}

From the above proof we see that the atoms in fact obey the BMO bound
$\| a_{I_T} \|_{BMO} \lesssim |I|^{-1/p}$.
One can improve this BMO control to $L^\infty$ control by repeating the argument in the John-Nirenberg inequality (Lemma \ref{jn}).  Namely, one locates the maximal sub-intervals where the averages of $a_{I_T}$ are large and separates off those trees, leaving behind a bounded atom.  One then repeats the process until only bounded atoms remain, in the spirit of Lemma \ref{tree-samp}.  We omit the details.  

Let $f$ be in BMO. From \eqref{bmo-def}, \eqref{w-size}, we have
$$ \| f \|_{BMO} = \sup_T \frac{1}{|I_T|^{1/2}} \| \Pi_T f \|_2.$$
Clearly it suffices to take suprema over complete trees $T$, thus by \eqref{piti}
$$ \| f \|_{BMO} = \sup_I \frac{1}{|I|^{1/2}} (\int_I |f - [f]_I|^2)^{1/2}.$$
By duality we thus have
\be{f-dual}
\| f \|_{BMO} = \sup_I \sup_{a \in S_0(I): \| a\|_2 = 1 } |I|^{-1/2} |\langle f, a \rangle|
\end{equation}
or equivalently
$$ \| f \|_{BMO} = \sup \{ |\langle f, a \rangle|: a \hbox{ is a } H^1 \hbox{ atom} \}.$$
Thus, as is well known, BMO is the dual of $H^1$.

\section{The Carleson embedding theorem and paraproducts}\label{embedding-sec}

We now give a slight variant of the above method, in which one selects trees using the averages $[f]_I$ instead of the sizes.  This type of argument is of course very old, and the arguments here are by no means new.  On the other hand, this type of tree selection method is a special case of the ``mean selection'' algorithm used (together with a size selection algorithm) in the proof of Carleson's theorem in \cite{laceyt-carleson}.  

We begin with

\begin{lemma}[Carleson embedding theorem]\label{embed}  Let $\P$ be a collection of lacunary tiles, $a: \P \to \R^+$ be a function, and $1 < p < \infty$.  Then we have
$$ \sum_{P \in \P} a(P) |[f]_{I_P}|^p \lesssim \| a \|_{\size^*(\P)} \|f\|_p^p$$
for all locally integrable functions $f$, with the implicit constants depending on $p$.
\end{lemma}

\begin{proof}
We apply Lemma \ref{mean-select} repeatedly, starting with a sufficiently large $n$ and decrementing $n$ repeatedly.  This gives us a partition
$\P = \bigcup_{n \in \Z} \bigcup_{T \in \T_n} T \cup \P_{-\infty}$
where the $\T_n$ are as in Lemma \ref{mean-select}, and 
$\| f \|_{\mean^*(\P_{-\infty})} = 0$.
The contribution of $\P_{-\infty}$ is zero, so it suffices to control
$$ \sum_n \sum_{T \in \T_n} \sum_{P \in T} a(P) |[f]_{I_P}|^p.$$
If $P \in T \in \T_n$, then
$|[f]_{I_P}| \leq \|f\|_{\mean(P)} \leq \|f\|_{\mean^*(T)} \lesssim 2^n$.
From this and \eqref{size-def}, \eqref{size-max-def}, \eqref{cheb} we may estimate the previous by
\bas &\lesssim \sum_n \sum_{T \in \T_n} \sum_{P \in T} a(P) 2^{np}\\
&\lesssim \sum_n \sum_{T \in \T_n} \|a\|_{\size^*(\P)} |I_T| 2^{np}\\
&\lesssim \sum_n \|a\|_{\size^*(\P)} 2^{np} 2^{-n} \int_{|f| \gtrsim 2^n} |f| \\
&\sim \| a\|_{\size^*(\P)} \int |f|^p
\end{align*}
as desired.
\end{proof}

We can apply this theorem to various linear and bilinear operators.  To do this we shall need some notation.

For any sequence $(a_P)_{P \in \P^+}$ of real numbers, we define the \emph{wavelet multiplier} $W^{-1} a_P W$ from $\S_0$ to $\S_0$ by
$$ W^{-1} a_P W f := \sum_{P \in \P^+} a_P Wf(P) \phi_P.$$
Wavelet multipliers are the discrete analogue of pseudo-differential operators, with $a_P$ being the discrete analogue of a symbol $a(x,\xi)$.  Observe that if $a_P$ is bounded, then $W^{-1} a_P W$ is bounded on $L^2$ and also bounded on BMO.  One can of course extend the domain $W^{-1} a_P W$ from $\S_0$ to $\S$, although some of the algebra properties are lost in doing so (since $W$ is not injective on $\S$).

Let $f$, $g$ be elements of $\S$.  We define the ``high-low'', ``low-high'', and ``high-high'' paraproducts\footnote{The continuous counterparts would be something like $\int (Q_t f) (P_t g) \frac{dt}{t}$, $\int (P_t f) (Q_t g) \frac{dt}{t}$, and $\int (Q_t f) (Q_t g) \frac{dt}{t}$, where $Q_t$ is as before and $P_t$ is a suitable approximation to the identity at width $1/t$, e.g. $P_t := e^{t^2 \Delta}$.  The precise definition of a paraproduct is not standardized, for instance $\pi_{hh}$ is not considered a paraproduct in some texts.}
\bas
\pi_{hl}(f, g) &:= \sum_{P \in \P^+}  Wf(P) [g]_P \phi_P\\
\pi_{lh}(f, g) &:= \sum_{P \in \P^+}  [f]_P Wg(P) \phi_P\\
\pi_{hh}(f, g) &:= \sum_{P \in \P^+}  Wf(P) Wg(P) \frac{\chi_{I_P}}{|I_P|}.
\end{align*}
These paraproducts have the symmetries
\be{permute}
\begin{split}
\int \pi_{hh}(f,g) h &= \int \pi_{hl}(g,h) f = \int \pi_{lh}(h,f) g\\
= \int \pi_{hh}(g,f) h &= \int \pi_{hl}(f,h) g = \int \pi_{lh}(h,g) f\\
&= \sum_{P \in \P^+} Wf(P) Wg(P) [h]_P
\end{split}
\end{equation}
and can be expressed in terms of the Littlewood-Paley  square function $S$:
$$ \pi_{hl}(f,g) = S^*( g Sf ); \quad \pi_{lh}(f,g) = S^*(f Sg); \quad \pi_{hh}(f,g) = Sf \cdot Sg.$$
When $f, g$ have mean zero (i.e. $f,g \in \S_0$), then the paraproducts decompose the pointwise product operator:
\be{point}
fg = \pi_{hl}(f,g) + \pi_{lh}(f,g) + \pi_{hh}(f,g).
\end{equation}
To see this, it suffices by bilinearity to reduce to the case when $f = \phi_P$ and $g = \phi_Q$ for some $P, Q \in \P^+$.  If $I_P$ and $I_Q$ are disjoint then both sides are zero.  Thus there are only three cases: $P >' Q$, $P <' Q$, and $P=Q$.  In these three cases the reader may easily verify that $fg$ is equal to the high-low, low-high, or high-high paraproduct of $f$ and $g$ respectively, and that the other two paraproducts vanish.

We observe that the high-low and low-high paraproducts can be written as wavelet multipliers:
\be{hllh-mult}
\pi_{hl}(f,g) = W^{-1} [g]_P Wf; \quad \pi_{lh}(f,g) = W^{-1} [f]_P Wg.
\end{equation}
The high-high paraproduct cannot be written in this way, but we have the useful relationship
\be{hh-mult}
\pi_{hh}(W^{-1} a_P Wf,g) = \pi_{hh}(f, W^{-1} a_P W g).
\end{equation}

One can also write paraproducts using both lacunary and non-lacunary tiles, for instance
\be{tril}
int \pi_{hh}(fg) h = \sum_{I \hbox{ dyadic}} |I|^{-1/2}
\langle f, \phi_{P^+(I)} \rangle
\langle g, \phi_{P^+(I)} \rangle
\langle h, \phi_{P^0(I)} \rangle.
\end{equation}
The bilinear Hilbert transform turns out to have a similar expansion, but with the sum ranging over a larger collection of triples of tiles than the ones for paraproducts (specifically, the tiles need not be lacunary or non-lacunary, and range over a three-parameter family rather than a two-parameter one).  See e.g. \cite{laceyt1}, \cite{laceyt2}, \cite{thiele-walsh}, \cite{cct}, \cite{thiele}, \cite{mtt:biest}.

From the Carleson embedding theorem we have paraproduct estimates:

\begin{corollary}[$L^2 \times BMO \to L^2$ paraproduct estimates]\label{paraproduct}  We have
$$
\| \pi_{hl}(f,g) \|_2\lesssim \| f\|_2 \|g\|_\infty$$
and
$$ \| \pi_{hh}(f,g) \|_2, \| \pi_{lh}(f,g) \|_2  
\lesssim \| f\|_2 \|g\|_{BMO}$$
for all $f, g \in \S$. 
\end{corollary}

In other words, paraproducts map $L^2 \times L^\infty$ to $L^2$ (just as the pointwise product does), and the $L^\infty$ factor can be relaxed to BMO as long as one only considers high frequencies of a BMO function.  Note that the low frequency portion of a BMO function is somewhat ill-defined since a BMO function might only be determined up to a constant.

\begin{proof} 
The first bound follows from \eqref{hllh-mult} since $[g]_P$ is bounded by $\|g\|_\infty$.  To prove the second bound, it suffices by \eqref{permute} to consider $\pi_{lh}$.  By orthogonality we have
$$ \| \pi_{lh}(f,g)\|_2 = (\sum_{P \in \P^+} |Wg(P)|^2 |[f]_{I_P}|^2)^{1/2}.$$
The claim now follows from Carleson embedding (Lemma \ref{embed}).
\end{proof}

\begin{lemma}[$BMO \times BMO \to BMO$ paraproduct estimate]\label{para-bmo}  We have
$$ \| \pi_{hh}(f, g) \|_{BMO} \lesssim \| f\|_{BMO} \|g\|_{BMO}$$
for all $f, g \in \S$.
\end{lemma}

For the other paraproducts $\pi_{hl}$, $\pi_{lh}$ one must place the ``low'' factor in $L^\infty$ rather than BMO, as in Lemma \ref{paraproduct}; this is again an easy consequence of \eqref{hllh-mult}.

\begin{proof}
By Lemma \ref{jn-char} with $p=1$ it suffices to show
$\int |\Pi_{\Tree(I)} \pi_{hh}(f,g)| \lesssim |I|$
for all dyadic intervals $I$.  

Fix $I$, and expand the left-hand side as
$$ \int |\sum_{P \in \P^+} Wf(P) Wg(P) \Pi_{\Tree(I)}(\frac{\chi_{I_P}}{|I_P|})|.$$
The summand vanishes unless $P \in \Tree(I)$.  Thus we can write the above as
$$ \int |\Pi_{\Tree(I)} \sum_{P \in \Tree(I)} Wf(P) Wg(P) \frac{\chi_{I_P}}{|I_P|}|.$$
Since $\Pi_{\Tree(I)}$ is bounded on $L^1$, we can bound this by
$$ \lesssim \int |\sum_{P \in \Tree(I)} Wf(P) Wg(P) \frac{\chi_{I_P}}{|I_P|}|.$$
Putting the absolute values inside and performing the integration, we can bound this by
$$ \lesssim \sum_{P \in \Tree(I)} |Wf(P)| |Wg(P)|.$$
The claim then follows from Cauchy-Schwarz and \eqref{bmo-def}.
\end{proof}

\subsection{Weak-type estimates}\label{weak-sec}

We now show how to use the above machinery to prove $L^{p,\infty}$ paraproduct estimates, where $L^{p,\infty}$ is the weak $L^p$ (quasi-)norm
$$ \| f \|_{L^{p,\infty}} := \sup_{\lambda > 0} \lambda | \{ x: |f(x)| \geq \lambda \} |^{1/p}.$$
We need the following basic characterization of weak $L^p$ for $0 < p < \infty$:

\begin{lemma}\label{weak-lp}  Let $0 < p < \infty$ and $A > 0$.  Then the following statements are equivalent up to constants:
\begin{enumerate}
\renewcommand{\labelenumi}{(\roman{enumi})}
\item $\| f \|_{p,\infty} \lesssim A$.
\item For every set $E$ with $0 < |E| < \infty$, there exists a subset $E' \subset E$ with $|E'| \sim |E|$ and $|\langle  f, \chi_{E'} \rangle| \lesssim A|E|^{1/p'}.$
\end{enumerate}
Here $p'$ is defined by $1/p' + 1/p = 1$ (note that $p'$ can be negative!).
\end{lemma}

\begin{proof}
To see that (i) implies (ii), set
$$ E' := E \backslash \{ x: |f(x)| \geq C A |E|^{-1/p} \}.$$
If $C$ is a sufficiently large constant, then (i) implies $|E'| \sim |E|$, and the claim follows.

To see that (ii) implies (i), let $\lambda > 0$ be arbitrary and set 
$E := \{ x: \Re(f(x)) > \lambda \}$.  Then by (ii) we have
$$ \lambda |E| \sim \lambda |E'| \lesssim A |E|^{1/p'},$$
and (i) easily follows (replacing $\Re$ by $-\Re$, $\Im$, $-\Im$ as necessary).
\end{proof}

When $p > 1$ we can always set $E' = E$, and the above lemma then reflects the duality between $L^{p,\infty}$ and $L^{p',1}$.  However for $p \leq 1$ the freedom to set $E'$ to be smaller than $E$ is necessary (since $f$ need not be
locally integrable).

A typical application of Lemma \ref{weak-lp} is

\begin{proposition}[$L^p \times L^q \to L^{r,\infty}$ paraproduct estimates]\label{paraproduct-2}  We have
$$ \| \pi_{hl}(f,g) \|_{r,\infty} \lesssim \| f\|_p \|g\|_q$$
whenever $1 < p,q < \infty$ and $1/p + 1/q = 1/r$.  
Similarly for $\pi_{lh}$, $\pi_{hh}$.
\end{proposition}

Note that $r$ can be less than 1.  One can strengthen the weak $L^r$ to strong $L^r$ by multilinear interpolation (see e.g. \cite{bergh:interp}, \cite{cct}, \cite{janson}).  The continuous version of these dyadic paraproduct estimates can be found in, e.g. \cite{coifmanm1}-\cite{coifmanm6}; the version for $r<1$ was first proven in \cite{grafakos:kalton} (with some special cases in \cite{calderon}, \cite{coifmanm1}).  It is possible to obtain the continuous estimates from the dyadic ones via averaging arguments, but we shall not do so here.  

\begin{proof}  We first consider $\pi_{hl}$. 
We may normalize $\|f\|_p = \|g\|_q = 1$; we may assume that $f$ and $g$ are dyadic test functions.  Let $E$ be a measurable set with $0 < |E| < \infty$.  We need to find a set $E' \subset E$ with $|E'| \sim |E|$ such that
\be{jim}
|\langle \pi_{hl}(f,g), \chi_{E'} \rangle|
\lesssim |E|^{1/r'}.
\end{equation}
By rescaling (using the hypothesis $1/p + 1/q = 1/r$) we may take $|E| \sim 1$.  

We choose $E'$ as
$$ E' := E \backslash \{ x: M|f|^p(x) + M|g|^q(x) \geq C \}.$$
If $C$ is large enough, then $|E| \sim |E'|$ by the Hardy-Littlewood maximal  inequality (see e.g. \cite{stein}).

We wish to show \eqref{jim} with $|E| \sim 1$.  By \eqref{permute} it suffices to show that
$|\sum_{P \in \P} Wf(P) [g]_{I_P} W\chi_{E'}(P)| \lesssim 1$
for all convex collections $\P$ of tiles.

We may remove all tiles in $\P$ for which $I_P \cap E' = \emptyset$, since $W\chi_{E'}$ vanishes on these tiles.  For any remaining tile $P$ we then have
\be{av}
\int_{I_P} |f|^p \lesssim |I_P|
\end{equation} 
by construction of $E'$.  Similarly, we have 
$$ \| g\|_{\mean^*(\P)} = \sup_{P \in \P} [|g|]_{I_P} \lesssim [|g|^q]_{I_P}^{1/q} \lesssim 1.$$
Thus we reduce to showing
\be{mash}
\sum_{P \in \P} |Wf(P)| |W\chi_{E'}(P)| \lesssim 1.
\end{equation}
From \eqref{av} and the $L^p$ boundedness of the Littlewood-Paley  square function (see e.g. \cite{stein}) we have
$$ \int_{I_P} |S\Pi_{\Tree(I_P)}f|^p \lesssim
\int_{I_P} |\Pi_{\Tree(I_P)} f|^p \lesssim \int_{I_P} |f|^p \lesssim |I_P|,$$
for all $P \in \P$, thus
$$ \int_{I_P} (\sum_{P \in \Tree(I_P)} |Wf(P)|^2 \frac{\chi_{I_P}(x)}{|I_P|})^{p/2}
\lesssim |I_P|.$$
Applying Chebyshev's inequality and Corollary \ref{good-lambda} we thus see that
$\| |Wf|^2 \|_{\size^*(\P)} \lesssim 1$.
Also, we have
$$ \| |W\chi_{E'}|^2 \|_{\size^*(\P)} \leq \| \chi_{E'} \|_{BMO}^2 \leq 
\| \chi_{E'} \|_\infty^2 \leq 1.$$

Thus we can find an $n = O(1)$ such that
$$ \| |Wf|^2 \|_{\size^*(\P_n)} \leq 2^{2n} \hbox{ and } \| |W\chi_{E'}|^2 \|_{\size^*(\P_n)} \leq 2^{2np/s},$$
where have set $\P_n := \P$, and $s>1$ is an exponent close to 1 to be chosen later.

By a finite  number of applications of Lemma \ref{tree-select} with $a := |Wf|^2$ or $a := |W \chi_{E'}|^2$ we may partition
$\P_n = \bigcup_{T \in \T_n} T \cup \P_{n-1}$
where $\P_{n-1}$ is a convex collection of tiles such that
$$ \| |Wf|^2 \|_{\size^*(\P_{n-1})} \leq 2^{2(n-1)} \hbox{ and } \| |W\chi_{E'}|^2 \|_{\size^*(\P_{n-1})} \leq 2^{2(n-1)p/s}$$
and $\T_n$ is a collection of convex trees with disjoint spatial supports such that either
$$ \int_{I_T} |\Pi_T f|^p \gtrsim 2^{np} |I_T|$$
or
$$ \int_{I_T} |\Pi_\CZ \chi_{E'}|^s \gtrsim 2^{np} |I_T|$$
for each $T \in \T_n$.  From \eqref{pit-max} we thus have
$$ \int (\tilde Mf)^p + \int (\tilde M \chi_{E'})^s
\gtrsim 2^{np} \sum_{T \in \T_n} |I_T|;$$
from our assumptions on $f$, $\chi_{E'}$ and the Hardy-Littlewood maximal  inequality (see e.g. \cite{stein}) we thus have
$\sum_{T \in \T_n} |I_T| \lesssim 2^{-np}$.
We now return to \eqref{mash}, and estimate the contribution of the trees in $\T_n$ by
$$ \sum_{T \in \T_n} \sum_{P \in T} |Wf(P)| |W\chi_{E'}(P)|.$$
We apply Cauchy-Schwarz followed by the size control on $Wf$ and $W\chi_{E'}$ we may bound this by
$$ \sum_{T \in \T_n} |I_T| \| |Wf|^2 \|_{\size(T)}^{1/2} \| |W\chi_{E'}|^2 \|_{\size(T)}^{1/2}
\lesssim \sum_{T \in \T_n} |I_T| 2^n 2^{np/s} \lesssim 2^{-np} 2^n 2^{np/s}.$$

We now turn to the contribution of the tiles in $\P_{n-1}$.  We may iterate the above procedure, decomposing $\P_{n-1}$ into $\T_{n-1}$ and $\P_{n-2}$, and continue in this fashion until we are left with a collection of tiles $\P_{-\infty}$ with size zero, which we can discard.  Summing up, we can thus control the left-hand side of \eqref{mash} by
$\lesssim \sum_{n \leq O(1)} 2^{-np} 2^n 2^{np/s}$.
If one chooses $s$ sufficiently close to 1, then this sum converges, and we are done.

A similar argument handles $\pi_{hl}$.  The remaining paraproduct $\pi_{hh}$ then follows from \eqref{point} and H\"older's inequality (which is still valid for $r < 1$).
\end{proof}

One can modify the above argument to obtain the corresponding estimate for the bilinear Hilbert transform (in the Walsh model, at least); see e.g. \cite{laceyt1}, \cite{thiele}, \cite{mtt:biest}.  A difficulty in that case is that the tiles are no longer lacunary, and one cannot guarantee the spatial disjointness of the trees $T$ in $\T_n$.  However one can still make the trees essentially disjoint in phase space, but then one can only use $L^2$ estimates to control $\sum_{T \in \T_n} |I_T|$ instead of $L^p$ estimates.  Because of this, the above strategy only seems to work for the bilinear Hilbert transform when $r > 2/3$; it appears that one needs very different techniques to handle the remaining case $1/2 < r \leq 2/3$.

\section{Calderon-Zygmund operators, and $\CZ(b)$ theorems}\label{cz-sec}

We now consider the theory of Calder\'on-Zygmund operators, and specifically the aspects of the theory related to $\CZ(1)$ and $\CZ(b)$ theorems.

For simplicity we shall restrict ourselves to model operators of the form
$$ \CZ f(x) := \int K(x,y) f(y)\ dy$$
where $K$ is a locally integrable function on $(x,y) \in [0, 2^M] \times [0,2^M]$ which obeys the kernel condition
\be{kernel}
|K(x,y)| \lesssim \frac{1}{|x-y|}
\end{equation}
and the \emph{perfect dyadic Calder\'on-Zygmund conditions}
\be{perfect}
|K(x,y) - K(x',y)| + |K(y,x) - K(y,x')| = 0
\end{equation}
whenever $x,x' \in I$ and $y \in J$ for some disjoint dyadic intervals $I$ and $J$.  Equivalently, $K$ is constant on all rectangles $\{ I \times J: I, J \hbox{ are siblings} \}$.  We also impose the technical truncation conditions\footnote{The continuous analogue of these conditions would be that $K(x-y)$ vanishes when $|x-y| \leq 2^{-M}$ or $|x-y| \geq 2^M$.  It is possible to formulate $T(b)$ theorems which do not require truncated operators, but this introduces some additional technicalities which are not relevant to the present discussion, and so we have chosen to ignore the issue for expository reasons.} that $K$ vanishes on all diagonal squares $I \times I$ where $|I| = 2^{-M}$, and also vanishes on the squares $[0,2^{M-1}] \times [2^{M-1}, 2^M]$ and
$[2^{M-1}, 2^M] \times [0, 2^{M-1}]$.

We refer to operators $\CZ$ of the above form as \emph{perfect dyadic Calder\'on-Zygmund operators}.  These operators can be thought of as the dyadic analogue of truncated Calder\'on-Zygmund operators, where the cancellation conditions are perfect.  (For ordinary Calder\'on-Zygmund operators one can bound the left-hand side of \eqref{perfect} by something like $O(|x-x'|/|x-y|^2)$).  

Let $\CZ$ be a perfect dyadic Calder\'on-Zygmund operator.  From \eqref{perfect} we observe that $\CZ$ maps $\S$ to $\S$, and furthermore if $I$ is any dyadic interval and $f \in \S_0(I)$, then $\CZ f$ is supported in $I$.  Similarly for $\CZ^*$.  We shall use this cancellation heavily in the sequel.

Examples of perfect dyadic Calder\'on-Zygmund operators include multiplier operators $W^{-1} a_P W$ and paraproducts $f \mapsto \pi(a,f)$.  It turns out that these are essentially the only such operators.  Indeed, we have the well-known splitting (see \cite{bcr}, \cite{meyer:ondelettes}) of a perfect dyadic Calder\'on-Zygmund operator into three parts: the diagonal part, the $\CZ(1)$ paraproduct, and the $\CZ^*(1)$ paraproduct:

\begin{lemma}\label{spl}  If $\CZ $ is a perfect Calder\'on-Zygmund operator, then we have 
\be{decomp}
\CZ f \equiv W^{-1} \langle \CZ \phi_P, \phi_P \rangle W f
+ \pi_{hl}(\CZ(1), f) + \pi_{hh}(\CZ^*(1), f)
\end{equation}
for all $f \in \S$, where $f \equiv g$ denotes the statement\footnote{This is equivalent to $f-g$ being constant on $[0,2^M]$.} that $Wf = Wg$.
\end{lemma}

\begin{proof}
We need to show that
\be{split}
\begin{split}
\langle \CZ f, \phi_Q \rangle = 
&\langle \CZ \phi_Q, \phi_Q \rangle Wf(Q) \\
&+ W(\CZ(1))(Q) [f]_Q \\
&+ 
\sum_{P \in \P^+} W(\CZ^*(1))(P) Wf(P) [\phi_Q]_P
\end{split}
\end{equation}
for all $f \in \S$ and $Q \in \P^+$.  

When $f = 1$ the identity is clear, so by subtracting off a multiple of 1 if necessary we may assume $f \in \S_0$, thus $f = \sum_P Wf(P) \phi_P$.  We can then decompose
$$ \langle \CZ f, \phi_Q \rangle = \sum_{P \in \P^+} Wf(P) \langle \CZ \phi_P, \phi_Q \rangle.$$
If $I_P$ and $I_Q$ are disjoint then $\langle \CZ \phi_P, \phi_Q \rangle = 0$ since $\CZ $ is a perfect dyadic Calder\'on-Zygmund operator.  Thus we may partition the sum into the portions $P=Q$, $Q <' P$, or $P <' Q$.

The diagonal term $P=Q$ is the first term in \eqref{split}.  Now consider the $Q <' P$ portion.  We write this as
$$ \langle \sum_{P \in \P^+: Q <' P} Wf(P) \phi_P, \CZ^* \phi_Q \rangle.$$
The function $\CZ^* \phi_Q$ is supported on $I_Q$, while 
$$ \sum_{P \in \P^+: Q <' P} Wf(P) \phi_P = f - 
\sum_{P \in \P^+: Q \not <' P} Wf(P) \phi_P$$
is constant on $I_Q$ and has the same mean as $f$ on $I_Q$, we thus have
$$ \langle \sum_{P \in \P^+: Q <' P} Wf(P) \phi_P, \CZ^* \phi_Q \rangle =
\langle [f]_Q, \CZ^* \phi_Q \rangle = \langle \CZ(1), \phi_Q \rangle [f]_Q$$
which is the second term in \eqref{split}.

Finally, we consider the $P <' Q$ term. Observe that the $P$ summation in \eqref{split} vanishes unless $P <' Q$.  It thus suffices to show that
$$ Wf(P) \langle \CZ \phi_P, \phi_Q \rangle = W(\CZ^*(1))(P) Wf(P) [\phi_Q]_P.$$
But this follows since $\CZ \phi_P$ is supported on $I_P$ and $\phi_Q$ is constant on $I_P$.
\end{proof}

\begin{corollary}[Dyadic global $\CZ(1)$ theorem]\label{t1}  Let $\CZ $ be a perfect dyadic Calder\'on-Zygmund operator such that 
$$ \| \CZ(1)\|_{BMO}, \|\CZ^*(1) \|_{BMO} \lesssim 1$$
and we have the weak boundedness property
\be{wbp-t1}
|\langle \CZ \phi_P, \phi_P \rangle| \lesssim 1 \hbox{ for all } P \in \P^+.
\end{equation}
Then $\CZ$ is bounded on $L^2$.
\end{corollary}

Note that the converse of this theorem is easy: if $\CZ$ is bounded on $L^2$, then we certainly have the weak boundedness property, and by \eqref{f-dual}
$$ \|\CZ(1)\|_{BMO} = \sup_I \sup_{a \in \S_0(I): \|a\|_2 = 1} |I|^{-1/2} |\langle \CZ(1),a \rangle| = \sup_I \sup_{a \in \S_0(I): \|a\|_2 = 1} |I|^{-1/2} |\langle \CZ(\chi_I), a \rangle| \lesssim 1$$
and similarly for $\CZ^*(1)$.  Indeed we observe that $\CZ$ and $\CZ^*$ must map $L^\infty$ to BMO.

\begin{proof}
From the conditions on the kernel $K$ and duality it suffices to show that
$|\langle \CZ f, g \rangle| \lesssim \|f\|_2 \|g\|_2$
for all $f, g \in \S$.  By splitting $[0,2^M]$ into two intervals it suffices to show this when $f \in \S([0,2^{M-1}])$ or $f \in S([2^{M-1}, 2^M])$.

Without loss of generality we may assume that $f \in \S([0, 2^{M-1}])$.
From \eqref{decomp}, Corollary \ref{paraproduct}, and the hypotheses on $\CZ$ we see that this estimate holds for all $g \in \S_0$. Also, from the truncation hypothesis the claim trivially holds for $g =  \chi_{[2^{M-1},2^M]}$.  Since $\S$ is spanned by $\S_0$ and $\chi_{[2^{M-1},2^M]}$, the claim follows. 
\end{proof}

We can rephrase Corollary \ref{t1} as an equivalent ``local'' version:

\begin{corollary}[Dyadic local $\CZ(1)$ theorem]\label{t1-local-cor} 
Let $\CZ$ be a perfect dyadic Calder\'on-Zygmund operator such that 
\be{t1-local}
\| \CZ(\chi_{I_P})\|_{L^1(I_P)}, \| \CZ^*(\chi_{I_P})\|_{L^1(I_P)} \lesssim |I_P| \hbox{ for all } P \in \P^+.
\end{equation}
Then $\CZ$ is bounded on $L^2$.
\end{corollary}

\begin{proof}
From \eqref{t1-local} we see that
$|\int_{I_P \times I_P} K(x,y)\ dx dy| \lesssim |I_P|$.
From \eqref{kernel} we thus have that
$$ |\int_{I_P \times I_P} \phi_P(x) K(x,y) \phi_P(y)\ dx dy| \lesssim 1,$$
or in other words that \eqref{wbp-t1} holds.  By Corollary \ref{t1} and symmetry it thus suffices to show that $\CZ(1) \in BMO$, or equivalently (by Corollary \ref{jn-char} and duality, cf. \eqref{f-dual}) that
$$ |\langle \CZ(1), h_I \rangle| \lesssim |I| \|h_I \|_\infty$$
whenever $I$ is a dyadic interval and $h_I \in \S_0(I)$.  But since $\CZ^*(h_I)$ is supported on $I$, we have $\langle \CZ(1), h_I \rangle = \langle \CZ \chi_I, h_I \rangle$,
and the claim follows from \eqref{t1-local}.
\end{proof}

Note that the converse of the above Corollary is immediate from H\"older's inequality.  One can also deduce Corollary \ref{t1} from Corollary \ref{t1-local-cor}, but we leave this to the reader.

We shall consider generalizations of the global and local $\CZ(1)$ theorem next, after some preliminaries on accretivity.

\subsection{Accretivity and one-sided $\CZ(b)$ theorems}

Let $b \in S$ be a complex-valued function, and $\P \subseteq \P^+$ be a collection of tiles. We say that $b$ is \emph{pseudo-accretive on $\P$} if
\be{weak-accrete}
|[b]_P| \gtrsim 1 \hbox{ for all } P \in \P.
\end{equation}
If we in addition have the property
\be{strong-accrete}
|[b]_{P_l}|, |[b]_{P_r}| \gtrsim 1
\end{equation}
for the two children\footnote{For this to be well defined, $\P$ cannot contain any tiles with minimal length $|I_P| = 2^{-M}$.} $P_l$, $P_r$ of tiles $P \in \P$, we say that $b$ is \emph{strongly pseudo-accretive} on $\P$.
Note that we are \emph{not} assuming any $L^\infty$ control on $b$ in the above definition. 

If $b$ is pseudo-accretive on the entire tile set $\P^+$, we simply say that $b$ is \emph{pseudo-accretive}.  Examples of pseudo-accretive functions include the \emph{accretive} functions, for which $\Re\ b(x) \gtrsim 1$ for all $x \in [0, 2^M]$

The $\CZ(1)$ theorem can now be generalized to ``one-sided $\CZ(b)$ theorems'' in which we control $\CZ(b)$ and $\CZ^*(1)$.  We give the dyadic version of an argument of Semmes \cite{semmes} (who considered the case $\CZ(b) \in BMO$, $\CZ^*(1)=0$):

\begin{theorem}[One-sided global $\CZ(b)$ theorem]\label{t1b}  Let $b$ be a pseudo-accretive function on $\P^+$ with $\| b\|_{BMO} \lesssim 1$.
Let $\CZ$ be a perfect dyadic Calder\'on-Zygmund operator obeying the weak boundedness property \eqref{wbp-t1} such that 
$$ \| \CZ(b)\|_{BMO}, \|\CZ^*(1) \|_{BMO} \lesssim 1.$$
Then $\CZ$ is bounded on $L^2$.
\end{theorem}

\begin{proof}
From the dyadic $\CZ(1)$ theorem (Corollary \ref{t1}) it suffices to show that $\| \CZ(1) \|_{BMO} \lesssim 1$, or in other words that $|W(\CZ(1))|^2$ has bounded maximal size (i.e. is a Carleson measure).

From \eqref{decomp} and \eqref{hllh-mult} 
\be{w-split}
\CZ(b) \equiv W^{-1} \langle \CZ \phi_P,\phi_P \rangle W b + W^{-1} [b]_P W \CZ(1) + \pi_{hh}(\CZ^*(1), b).
\end{equation}
From \eqref{weak-accrete} we may thus solve for $\CZ(1)$:
\be{t1-form}
\CZ(1) \equiv W^{-1} [b]_P^{-1} W [\CZ(b) - W^{-1} \langle \CZ \phi_P,\phi_P \rangle W b - \pi_{hh}(\CZ^*(1), b)]
\end{equation}
and the claim follows from \eqref{weak-accrete}, \eqref{wbp-t1}, the hypotheses $\|\CZ(b)\|_{BMO}, \|b\|_{BMO} \lesssim 1$, and Lemma \ref{paraproduct}.
\end{proof}

Note that one only needs $b$ in BMO in the above argument instead of the more usual $L^\infty$.

One drawback to the above theorem is that it requires the function $b$ to be pseudo-accretive.  Fortunately, it is not too difficult to construct pseudo-accretive functions.  The following basic lemma says that if a function has large mean, then it is pseudo-accretive on a non-trivial set of tiles.  (Equivalently, if a function has small mean on too many small tiles, then it must have small mean globally).

\begin{lemma}\label{accrete-lemma}  Let $T_0 \subseteq \P^+$ be a convex tree, and let $b$ be a function such that
\be{2-bound}
\| \Pi_{T_0} b \|_2 \leq C_0 |I_{T_0}|^{1/2}
\end{equation}
and
$ |[b]_{I_{T_0}}| \geq \delta$
for some $C_0, \delta > 0$.  Then there exists $0 < \eps \ll 1$ depending only on $C_0$ and $\delta$ and a family $\T$ of disjoint convex sub-trees of $T_0$ whose tops form a $(1-\eps)$-packing of $T_0$, and such that
$|[b]_P| > \eps$
for all $P \in T_0 \backslash \bigcup_{T \in \T} T$.  Furthermore we have
$|[b]_{P_T}| \leq \eps \hbox{ for all } T \in \T$.
\end{lemma}

\begin{proof}
Let $\P$ denote those tiles in $T_0$ for which $|[b]_P| \leq \eps$, and which are maximal with respect to the ordering $<'$.  Clearly the tiles in $\P$ have disjoint spatial intervals and obey
$$ |\int_{I_P} b| < \eps |I_P|.$$
To prove the lemma it will suffice to show the $(1-\eps)$-packing property. 
Suppose for contradiction that
\be{contra}
|\bigcup_{P \in \P} I_P| > (1-\eps) |I_{T_0}|.
\end{equation}
Using the identity
$\int_{I_P} b = \int_{I_P} ([b]_{I_{T_0}} + \Pi_{T_0} b)$
and summing in $P$ we obtain
$$ |\int_{\bigcup_{P \in \P} I_P} [b]_{I_{T_0}} + \Pi_{T_0} b| \leq \eps |\bigcup_{P \in \P} I_P| \leq \eps |I_{T_0}|;$$
if $\eps$ is sufficiently small with respect to $\delta$, we thus see that
$$|\int_{\bigcup_{P \in \P} I_P} \Pi_{T_0} b| \geq \frac{\delta}{2} |I_{T_0}|.$$
Since $\Pi_{T_0} b$ has mean zero on $I_{T_0}$, we thus have
$$|\int_{I_{T_0} \backslash \bigcup_{P \in \P} I_P} \Pi_{T_0} b| \geq \frac{\delta}{2} |I_{T_0}|.$$
But this will contradict \eqref{contra} and \eqref{2-bound} by Cauchy-Schwarz, if $\eps$ is sufficiently small.
\end{proof}

This lemma combines nicely with Lemma \ref{tree-samp}.  Together, these lemmata heuristically allow us to treat ``large mean'' as being equivalent to ``pseudo-accretive'', at least for the purposes of placing something in BMO.
As an application we now use this lemma to give a localized version of Theorem \ref{t1b}.  Similar results\footnote{More precisely, the solution to the Kato problem requires a matrix-valued analogue of this theorem in which the $b_P$ are matrix valued, and $\CZ$ maps matrices to vectors.  The argument then requires an additional subtlety, namely a preliminary partition of the tile set $\P^+$ into $O(1)$ pieces, where on each piece the vector-valued coefficients $W(\CZ({\bf 1}))(P)$ lie in a narrow conical region.  See \cite{kato-2}.} have been used to solve the Kato problem in higher dimensions (see e.g. \cite{kato-2}):

\begin{theorem}[Local one-sided $\CZ(b)$ theorem]\label{local-tb1}  Let $\CZ$ be a perfect Calder\'on-Zygmund operator obeying $\| \CZ^*(1) \|_{BMO} \lesssim 1$ and \eqref{wbp-t1}.
Suppose also that for every $P \in \P^+$ there exists a function $b_P \in \S(I_P)$ with the normalization condition
$$
[b_P]_P = 1$$
and the $L^2$ bounds
$$
\int_{I_P} |b_P|^2 + |\CZ b_P|^2 \lesssim |I_P|
$$
for all $P \in \P^+$.
Then $\CZ$ is bounded on $L^2$.
\end{theorem}

The weak boundedness condition \eqref{wbp-t1} can actually be removed; see the remarks after Theorem \ref{ctb}.  Informally, this theorem asserts that to prove the $L^2$ boundedness of an operator $\CZ$ it actually suffices to establish boundedness for a single function $b_P$ for each interval $I_P$, provided that $b_P$ is not degenerate (in the sense that its mean is large) and provided that $T^*(1)$ is under control.  In the next section we shall remove the condition on $T^*(1)$, obtaining a ``two-sided'' version of this theorem.

\begin{proof}
Again it suffices to show that $\CZ(1)$ is in BMO.  By Lemma \ref{tree-samp} it suffices to show that for every complete tree $T$ we have
\be{weak-targ}
\sum_{P \in T \backslash \bigcup_{T' \in \T} T'} |W(\CZ(1))(P)|^2 \lesssim |I_T|
\end{equation}
for some collection $\T$ of disjoint convex trees in $T$ whose tops form a $(1-\eps)$-packing of $T$ for some $\eps > 0$.

Fix $T$.  By Lemma \ref{accrete-lemma} we can indeed find such a collection $\T$ with the additional property that $b$ is pseudo-accretive on $T \backslash \bigcup_{T' \in \T} T'$.  The claim \eqref{weak-targ} then follows from the argument used to prove Corollary \ref{t1b}.
\end{proof}

The above arguments do not extend well to two-sided situations in which one controls $\CZ(b_1)$ and $\CZ^*(b_2)$ (unless one of $b_1$, $b_2$ is close to a constant, e.g. in BMO norm).  In order to handle the general case we need adapted Haar bases, to which we now turn.
 
\subsection{Adapted Haar bases, and two-sided $\CZ(b)$ theorems}

Let $\P$ be a collection of tiles, and let $b$ be a function which is strongly pseudo-accretive on $\P$.  For each $P \in \P$, we define the adapted Haar wavelet $\phi^b_P$ (introduced in \cite{cjs}; see also \cite{at}) by
\be{phib-ident}
\phi^b_P := 
|I_P|^{-1/2} \frac{[b]_{P_r}}{[b]_{P}} \chi_{I_{P_l}} -
|I_P|^{-1/2} \frac{[b]_{P_l}}{[b]_{P}} \chi_{I_{P_r}}.
\end{equation}
Observe that this collapses to $\phi_P$ if $b$ is constant on $I_P$.  For non-constant $b$, $\phi^b_P$ is no longer mean zero, but one can easily verify that $\phi^b_P$ still obeys the weighted mean zero condition
\be{mean}
\int b \phi^b_P = 0.
\end{equation}
As a consequence we have the orthogonality property
\be{ortho}
\int \phi^b_P b \phi^b_Q = 0 \hbox{ for all distinct } P, Q \in \P.
\end{equation}
From \eqref{phib-ident} we see that
\be{phi-calc}
\int \phi^b_P b \phi^b_P 
= \frac{[b]_{P_r} [b]_{P_l}}{[b]_P}
= \frac{2}{[b]_{P_l}^{-1} + [b]_{P_r}^{-1}}.
\end{equation}
In particular, from the strong pseudo-accretivity condition \eqref{strong-accrete} we have the bound 
\be{norm}
|\int \phi^b_P b \phi^b_P| \gtrsim 1.
\end{equation}
It is interesting that this bound uses only the strong pseudo-accretivity of $b$, and in particular does not require $L^\infty$ control on $b$.

Define the dual adapted Haar wavelet $\psi^b_P$ by
$$ \psi^b_P := \frac{\phi^b_P b}{\int \phi^b_P b \phi^b_P}.$$
By \eqref{ortho}, \eqref{norm} we thus have that
$
\langle \psi^b_P, \phi^b_Q \rangle = \delta_{PQ}
$
where $\delta$ is the Kronecker delta.  In particular we have the representation formula
\be{g-expand}
f = \sum_{P \in \P} W_b f(P) \psi^b_P
\end{equation}
whenever $g$ is in the span of $\{ \psi^b_P: P \in \P \}$, where the adapted wavelet coefficients $W_b f(P)$ are defined by
$$ W_b f := \langle f, \phi^b_P \rangle.$$

We have the following basic orthogonality property:

\begin{lemma}\label{ortho-lemma}  Let $T$ be a convex tree, and let $b$ be a function which is pseudo-accretive on $T$ and obeys the mean bound 
\be{b-bmo}
\| |b|^2 \|_{\mean^*(T)} \lesssim 1.
\end{equation}
Then for any function\footnote{It can easily be seen, by aid of \eqref{g-expand}, that the estimate \eqref{f-ortho} can be reversed for all $f$ in the span of the $\phi^b_P$, but we will not use this.} $f \in S$ we have 
\be{f-ortho}
(\sum_{P \in T} |W_b f(P)|^2)^{1/2} \lesssim \|f\|_2.
\end{equation}
In fact, the more general estimate 
\be{f-var}
(\sum_{P \in T} |W_b (b' f) (P)|^2)^{1/2} \lesssim \|f\|_2 \| |b'|^2 \|_{\mean^*(T)}^{1/2}
\end{equation}
holds for any $f, b' \in S$.
\end{lemma}

\begin{proof}
From \eqref{b-bmo} and \eqref{w-size} we observe that 
\be{b-bmo-2}
\| |Wb|^2\|_{\size^*(T)} \lesssim 1.
\end{equation}

We first prove \eqref{f-ortho}.  From \eqref{phib-ident} and the identity
$|I_P|^{-1/2} Wb(P) = [b]_{P_l} - [b]_P = [b]_P - [b]_{P_r}$, we obtain the identity 
$$
W^b f(P) = W f(P) - \frac{Wb(P)}{[b]_P} [f]_P.
$$
If we replace $W_b$ by $W_b$ then \eqref{f-ortho} follows from Bessel's inequality and the orthonormality of the Haar wavelets $\phi_P$.  By the previous identity and the triangle inequality it thus suffices to show
$$ (\sum_{P \in T} |\frac{Wb(P)}{[b]_P} [f]_P|^2)^{1/2} \lesssim \|f\|_2.$$
We may discard $[b]_P$ by pseudo-accretivity \eqref{weak-accrete}.  The claim then follows from Carleson embedding (Lemma \ref{embed}) and \eqref{b-bmo-2}.

Now we prove \eqref{f-var}.  Let $\Q$ denote the collection of tiles in $\P^+$ which are children of tiles in $T$, but are not in $T$ itself.  In order to ensure that the intervals $\{ I_Q: Q \in \Q \}$ partition $I_T$ we will allow the tiles $Q$ to have spatial intervals $|I_Q| = 2^{-M-1}$; the partition property then follows from the convexity of $T$.  

The function
$ b'f - \sum_{Q \in \Q} [b'f]_Q \chi_{I_Q}$
has mean zero on every interval $I_Q$, and is thus orthogonal to $\phi^b_P$ for every $P \in T$.  We may thus freely replace $b'f$ by the averaged function $\sum_{Q \in \Q} [b'f]_Q \chi_{I_Q}$ in \eqref{f-var}.  By \eqref{f-ortho} it thus suffices to show that
$$ \| \sum_{Q \in \Q} [b'f]_Q \chi_{I_Q}\|_2^2 \lesssim \| f\|_2^2 \| |b'|^2 \|_{\mean^*(T)}.$$
But from Cauchy-Schwarz we have
$$ |[b'f]_Q|^2 |I_Q| \lesssim \| f\|_{L^2(I_Q)}^2 \| |b'|^2 \|_{\mean(2Q)}
\leq \| f\|_{L^2(I_Q)}^2 \| |b'|^2 \|_{\mean^*(T)},$$
and the claim follows by summing in $Q$.
\end{proof}

We can now give our main result, namely a dyadic local $\CZ(b)$ theorem. 

\begin{theorem}[Dyadic local $\CZ(b)$ theorem]\label{ctb}  Let $\CZ$ be a perfect Calder\'on-Zygmund operator, and suppose that for each $P \in \P^+$ we can find functions $b^1_P$, $b^2_P$ in $\S(I_P)$ obeying the normalization
\be{normal}
[b^1_P]_P = [b^2_P]_P = 1
\end{equation}
and the bounds
\be{bounds}
\int_{I_P} |b^1_P|^2 + |\CZ b^1_P|^2 + |b^2_P|^2 + |\CZ^*b^2_P|^2 \lesssim |I_P|.
\end{equation}
Then $\CZ$ is bounded on $L^2$.
\end{theorem}

This theorem is a stronger version of the local $\CZ(b)$ theorem in 
\cite{christ-tb} (but for the dyadic setting with perfect cancellation), which required $L^\infty$ control in \eqref{bounds} instead of $L^2$ control. (This was generalized to BMO control and to non-doubling situations in \cite{ntv:tb-local}).  Also it required the global $\CZ(b)$ theorem of David, Journ\'e and Semmes \cite{djs} (which we instead deduce as a corollary of Theorem \ref{ctb}). We make some further remarks after the proof of the theorem.

\begin{proof}
This proof is somewhat lengthy and so we split the argument into several stages.

{\bf Step 0.  Preliminary estimates.}

We begin with a basic lemma which already shows the importance of the normalization \eqref{normal}.

\begin{lemma}[$b^1_P$ spans $S(I_P)/S_0(I_P)$]\label{split-lemma}  For any tile $P \in \P^+$ and any $f \in S(I_P)$, we have
$$ \| f \|_{L^2(I_P)} \lesssim \| f - [f]_P \|_{L^2(I_P)} + |I_P|^{-1/2} |\langle f, b^1_P \rangle|.$$
Similarly for $b^2_P$.
\end{lemma}

\begin{proof}
Let $h$ be an arbitrary element of $S(I_P)$ with $\|h\|_2 = 1$.  Then
$$ \langle f, h \rangle = \langle f, h - [h]_P b^1_P \rangle + [h]_P \langle f, b^1_P \rangle = \langle f - [f]_P, h - [h]_P b^1_P \rangle + [h]_P \langle f, b^1_P \rangle.$$
By Cauchy-Schwarz and \eqref{bounds} we thus have
$$ |\langle f, h \rangle| \lesssim \| f - [f]_P \|_{L^2(I_P)} + |I_P|^{-1/2} |\langle f, b^1_P \rangle|.$$
Taking suprema over all $h$, the claim follows.
\end{proof}

A useful application of the above lemma is the following convenient truncation property of the $b^1_P$ and $b^2_P$ (already observed in \cite{christ-tb}):

\begin{corollary}\label{trunc}  Let $P, Q$ be lacunary tiles with $Q \leq' P$.  If we have the estimate
\be{hyp}
\int_{I_Q} |\CZ b^1_P|^2 + |b^1_P|^2 \lesssim K |I_Q|
\end{equation}
for some $K \gtrsim 1$, then we have
$$ \int_{2I_Q} |\CZ(b^1_P \chi_{I_Q})|^2 \lesssim K |I_Q|.$$
Similarly for $b^2_P$ (but with $\CZ$ replaced by $\CZ^*$).
\end{corollary}

\begin{proof}
By \eqref{kernel} and \eqref{hyp} the portion of the integral on $2I_Q \backslash I_Q$ is acceptable, so it suffices to bound the integral on $I_Q$.
From Cauchy-Schwarz, \eqref{bounds} and \eqref{hyp} we have
$$ |\langle \CZ(b^1_P \chi_{I_Q}), b^2_Q \rangle| = |\langle b^1_P \chi_{I_Q}, \CZ^* b^2_Q \rangle| \leq \| b^1_P \|_{L^2(I_Q)} \| \CZ^* b^2_Q \|_{L^2(I_Q)} \lesssim K^{1/2} |I_Q|.$$
By Lemma \ref{split-lemma} it thus suffices to show that
$$ \| \CZ(b^1_P \chi_{I_Q}) - [\CZ(b^1_P \chi_{I_Q})]_Q \|_{L^2(I_Q)} \lesssim K^{1/2} |I_Q|^{1/2}.$$
Now observe that for every $h \in S_0(I_Q)$ we have
$$ \langle \CZ(b^1_P \chi_{I_Q}), h \rangle = \langle b^1_P \chi_{I_Q}, \CZ^* h \rangle = \langle b^1_P, \CZ^* h \rangle = \langle \CZ b^1_P, h \rangle.$$
By duality this implies that
$$ \CZ(b^1_P \chi_{I_Q}) - [\CZ(b^1_P \chi_{I_Q})]_Q = \CZ(b^1_P) - [\CZ(b^1_P)]_Q$$
on $I_Q$.  The claim then follows from \eqref{hyp}.
\end{proof}

This Corollary will be useful in estimating the operator $\CZ$ when acting on objects such as $\psi^{b^1_P}_Q$ which can be expressed as linear combinations of truncated versions of $b^1_P$.  Similarly when estimating $\CZ^*$ on objects such as $\psi^{b^2_P}_Q$.

{\bf Step 1.  Overview of main argument.}

We now begin the main argument.  Let $A$ be the best constant such that
$$ \| \CZ^* \chi_{I_P} \|_{L^1(I_P)} \leq A |I_P|$$
for all tiles $P \in \P^+$.  We claim that $A = O(1)$; from this and the corresponding claim for $\CZ \chi_{I_P}$ (which is of course symmetric) the theorem will follow from the local $\CZ(1)$ theorem (Corollary \ref{t1-local-cor}).

In fact we will show 
\be{dual-targ}
|\int_{I_P} \CZ f| \leq ((1-\eps)A + O(1)) |I_P| \|f\|_\infty
\end{equation}
for all tiles $P \in \P^+$ and $f \in S(I_P)$, and some $0 < \eps \ll 1$ depending only on the implicit constant in \eqref{bounds}.  By duality this implies that $A \leq (1-\eps)A + O(1)$, which will prove the desired bound on $A$.

Fix $P$, $f$.  We shall prove the estimate \eqref{dual-targ} in three steps.  Firstly (in Step 2), we decompose $f$ and reduce matters to proving a Carleson measure type estimate on the wavelet coefficients $|\langle \CZ^* \chi_{I_P}, \psi^{b^1_P}_Q \rangle|^2$; this argument shall use stopping-time arguments (which we encapsulate as Lemma \ref{subtree}) based on $b^1_P$ but not on $b^2_P$.  Then (in Step 3), we decompose $\chi_{I_P}$ and use stopping time arguments (again using Lemma \ref{subtree}) based on $b^2_P$ but not on $b^1_P$.  It will be important not to try to handle $b^1_P$ and $b^2_P$ at the same time as we will lose the crucial $(1-\eps)$ packing property of the trees left out by the stopping time algorithm if we do so.  

The purpose of these stopping arguments is to impose some pseudo-accretivity and other regularity properties on the $b^1_P$ and $b^2_P$.  Once we have enough regularity properties, we can then (in Step 4) do an elementary computation to estimate the wavelet coefficients $|\langle \CZ^* \chi_{I_P}, \psi^{b^1_P}_Q \rangle|$ pointwise by the quantities which we know to be controlled by hypothesis (see \eqref{bounds} below).

{\bf Step 2.  Pruning the bad tiles of $b^1_P$.}

We now begin the first of the three steps outlined above. 
We would like to break up $f$ into linear combinations of the wavelets $\psi^{b^1_P}_Q$, but we cannot do this for all $Q$ because we do not control the strong pseudo-accretivity of $b^1_P$.  However, by using Lemma \ref{accrete-lemma} and some other selection algorithms we can find a large subtree of $\Tree(P)$ for which we can decompose $f$ as desired, modulo acceptable errors:

\begin{lemma}\label{subtree}  Let $P \in \P^+$ be a tile.  Then we can partition
$$ \Tree(P) = T_1 \cup \P_{buffer} \cup \bigcup_{T' \in \T} T'$$
where 
\begin{itemize}
\item $\T$ is a collection of disjoint complete trees in $\Tree(P)$ whose tops form a $(1-\eps)$-packing of $\Tree(P)$ for some $0 < \eps \ll 1$ (depending only on the implicit constant in \eqref{bounds});
\item $T_1$ is a tree with top $P$ such that $b^1_P$ is strongly pseudo-accretive on $T_1$ (with constants perhaps depending on $\eps$);
\item $\P_{buffer}$ is a $2$-packing of $\Tree(P)$, $T_1 \cup \P_{buffer}$ is convex, $b^1_P$ is pseudo-accretive on $T_1 \cup \P_{buffer}$ and we have the mean bounds
\be{mean-bounds}
\| |b^1_P|^2 + |\CZ b^1_P|^2 \|_{\mean^*(T_1 \cup \P_{buffer})} \lesssim 1
\end{equation}
(with the implicit constant depending on $\eps$).
\item We have the decomposition
\be{f-form}
f = [f]_P b^1_P + \sum_{Q \in T_1} W_{b^1_P} f(Q) \psi^{b^1_P}_Q + \sum_{T' \in \T} (f \chi_{I_{T'}} - [f]_{P_{T'}} b^1_{P_{T'}})
+ \sum_{Q \in \P_{buffer}} \varphi_Q
\end{equation}
whenever $f \in S(I_P)$, where the ``buffer functions'' $\varphi_Q$ are  supported on $I_Q$, have mean zero, and take the form
$$ \varphi_Q = a_Q b^1_P \chi_{I_{Q_l}} + a'_Q b^1_P \chi_{I_{Q_r}} + a''_Q b^1_{Q_l} + a'''_Q b^1_{Q_r}$$
where the co-efficients $a_Q$, $a'_Q$, $a''_Q$, $a'''_Q$ depend on $f$ and the $b^1_P$ and (when $|I_Q| \neq 2^{-M}$) obey the bounds
\be{aaaa}
|a_Q| + |a'_Q| + |a''_Q| + |a'''_Q| \lesssim \|f\|_\infty.
\end{equation}
\end{itemize}
A similar statement holds with $b^1_P$ and $\CZ b^1_P$ replaced by $b^2_P$ and $\CZ^* b^2_P$ (but the sets $T_1$, $\P_{buffer}$ and $\T$ are different then).
\end{lemma}

The tree $T_1$ represents the ``good'' portion of the tree $\Tree(P)$, in which $b^1_P$ is neither too large nor too small (so in particular the $W^{b^1_P}$
wavelet system is well-behaved on $T_1$).  The buffer tiles $P_{buffer}$ 
are those tiles immediately above $T_1$ 
(and are thus slightly less ``good''), while the remaining trees 
$\T$ have no good properties at all, except that they only occupy at 
most $(1-\eps)$ of the tree $\Tree(P)$.  This decomposition shares
many features in common with Lemma \ref{slice} (for instance, the
trees $\T$ are formed from those intervals where $b$ is too ``heavy'' or
too ``light'').

In terms of the phase plane (but adapted to the $W_{b^1_P}$ wavelet system instead of the Haar wavelet system), one can interpret the right-hand side of \eqref{f-form} as follows.  The first term corresponds to the region of phase space below the tree $T_1$.  The second term corresponds to $T_1$ itself.  The third term corresponds to the region above $T_1 \cup \P_{buffer}$, while the last term is an error term corresponding to the region $\P_{buffer}$.  In the model case $b^1_P = \chi_{I_P}$, \eqref{f-form} simplifies to
$$ f = [f]_P \chi_{I_P} + \Pi_{T_1} f + \sum_{T' \in \T} \Pi_{T'} f + \sum_{Q \in \P_{buffer}} Wf(Q) \phi_Q.$$ 

\begin{proof}
We begin by applying Lemma \ref{accrete-lemma} to $\Tree(P)$ to find a preliminary collection $\T_0$ of disjoint convex trees in $\Tree(P)$ such that the tops of $\T_0$ are a $(1-2\eps)$-packing, and such that $b^1_P$ is pseudo-accretive (with constants depending on $\eps$) on the tree
$$ T_2 := \Tree(P) \backslash \bigcup_{T' \in \T_0} T'.$$
However we do not yet have \eqref{mean-bounds}.  To obtain these bounds we let $\Q$ denote the set of all tiles $Q \in T_2 $ for which
$$ \| |b^1_P|^2 + |\CZ b^1_P|^2 \|_{\mean(Q)} \geq C/\eps$$
and which are maximal with respect to $\leq'$.  If the constant $C$ is chosen large enough, then $\Q$ is a $\eps$-packing of $\Tree(P)$.  Thus if we define
$$ \T := \T_0 \cup \bigcup_{Q \in \Q} (\Tree(Q) \cap T_2)$$
then we see that \eqref{mean-bounds} holds on the tree
$$ T_3 := \Tree(P) \backslash \bigcup_{T' \in \T} T',$$
while the tops of $\T$ are still a $(1-\eps)$-packing of $\Tree(P)$.

We now perform one minor modification to $\T$ to make $\T$ sibling-free.  If $\T$ contains two trees whose tops $P_{T'}$, $P_{T''}$ are siblings, we can concatenate these trees and add a new tile $2P_{T'} = 2P_{T''}$ to join these trees to a larger tree without affecting the $(1-\eps)$-packing nature of the tree tops.  Repeating this process as often as necessary (it must terminate since $\Tree(P)$ only has a finite number of tiles) we can make $\T$ sibling-free.

For similar reasons we may assume that the trees\footnote{Alternatively, we could avoid these modifications by combining the stopping time argument here with the one in Lemma \ref{accrete-lemma}.} in $\T$ are complete, since we can always replace an incomplete tree by the completion of that tree, absorbing any sub-trees that were also in $\T$ if necessary.

We now define\footnote{The algorithm here is extremely similar to the one used to prove Theorem \ref{tree-slice}.  Indeed, one can even re-use Figure \ref{extrap-fig}.  The trees $\T_0$ are the ``light'' trees where $b^1_P$ has too small a mean; the tiles $\Q$ correspond to the circled ``heavy'' tiles, where $b^1_P$ or $\CZ b^1_P$ has too large an $L^2$ norm.  The tiles $\P_{buffer}$ are thus the buffer tiles, which are the ones just below the heavy or light tiles, as well as the tiles at the very finest scale.} $\P_{buffer}$ to be the set of tiles $Q$ in $T_3$ such that one or both\footnote{Of course, because we made $\T$ sibling-free, the only way both the children of $Q$ fail to be in $T_3$ is if $Q$ is at the finest scale, i.e. if $|I_Q| = 2^{-M}$.} of the children $Q_l$, $Q_r$ of $Q$ are not in $T_3$.   Since $T_3$ is a convex tree, the children of $Q$ who are not in $T_3$ must have disjoint spatial supports as $Q$ varies in $\P_{buffer}$.  This implies that $\P_{buffer}$ is a 2-packing.  We now set $T_1 := T_3 \backslash \P_{buffer}$.  Note that all children of tiles in $T_1$ lie in $T_1 \cup \P_{buffer}$ so that $b^1_P$ is strongly pseudo-accretive on $T_1$, but is merely pseudo-accretive on $T_1 \cup \P_{buffer}$.

The only property left to verify is the decomposition \eqref{f-form}.  If $f$ is a constant multiple of $b^1_P$ then only the first term is non-zero (thanks to \eqref{mean}) and the claim is easily verified.  By subtracting off a constant multiple we may thus assume that $f$ has mean zero on $I_P$.

It will suffice to prove the identity assuming that
$$ [b^1_P]_Q, [b^1_P]_{Q_l}, [b^1_P]_{Q_r} \neq 0 \hbox{ for all } Q \in \Tree(P),$$
since the general case then follows by an obvious limiting argument (the bounds \eqref{aaaa} will not depend quantitatively on the above condition).  In this case \eqref{g-expand} applies\footnote{One can easily verify (either by a dimension counting argument, or by inductively working from the finest scale upwards) that the wavelets $\psi^{b^1_P}_Q$ for $Q \in \Tree(P)$ span $S_0(P)$.}.  Comparing this with \eqref{f-form} and using the mean zero condition, we reduce to showing that
$$ 
\sum_{Q \in \P_{buffer}} W_{b^1_P} f(Q) \psi^{b^1_P}_Q + \sum_{T' \in \T} \sum_{Q' \in T'} W_{b^1_P} f(Q') \psi^{b^1_P}_{Q'} 
= \sum_{T' \in \T} (f \chi_{I_{T'}} - [f]_{P_{T'}} b^1_{P_{T'}}) + \sum_{Q \in \P_{buffer}} \varphi_Q$$
for suitable $\varphi_Q$.

Let $Q \in \P_{buffer}$.  First suppose that neither child of $Q$ is a top of a tree in $\T$ (since $\Tree(P)$ is complete, this can only happen when $|I_Q| = 2^{-M}$).  In this case we simply set $\varphi_Q := W_{b^1_P} f(Q) \psi^{b^1_P}_Q$.

Now suppose that one child of $Q$ is a top of a tree $T'$ in $\T$; without loss of generality we assume $Q_l$ is such a top.  Since $\T$ is sibling-free, $Q_r$ is not a top and must therefore lie in $T_1 \cup \P_{buffer}$.  In particular we have the lower bounds
\be{qqr}
|[b^1_P]_{Q_r}|, |[b^1_P]_Q| \gtrsim 1.
\end{equation}
We do not have good lower bounds on $|[b^1_P]_{Q_l}|$, but fortunately we can 
take advantage of some ``wiggle room'' in the buffer, and 
avoid using this in our computations by exploiting the identity 
$$ W_{b^1_P} f(Q) \psi^{b^1_P}_Q + \sum_{Q' \in T'} W_{b^1_P} f(Q') \psi^{b^1_P}_{Q'} = \sum_{Q' \in \Tree(Q)} W_{b^1_P} f(Q') \psi^{b^1_P}_{Q'}
- \sum_{Q' \in \Tree(Q_r)} W_{b^1_P} f(Q') \psi^{b^1_P}_{Q'}.$$

The function $\sum_{Q' \in \Tree(Q)} W_{b^1_P} f(Q') \psi^{b^1_P}_{Q'}$ clearly is supported on $I_Q$ and has mean zero, while
$$ f - \sum_{Q' \in \Tree(Q)} W_{b^1_P} f(Q') \psi^{b^1_P}_{Q'}  = \sum_{Q' \in \Tree(P) \backslash \Tree(Q)} W_{b^1_P} f(Q') \psi^{b^1_P}_{Q'} $$
is a constant multiple of $b^1_P$ on $I_{T'}$.  Thus we have
$$ \sum_{Q' \in \Tree(Q)} W_{b^1_P} f(Q') \psi^{b^1_P}_{Q'} = f \chi_{I_Q} - \frac{[f]_Q}{[b^1_P]_{Q}} b^1_P \chi_{I_Q}.$$
Similarly we have
$$ \sum_{Q' \in \Tree(Q_r)} W_{b^1_P} f(Q') \psi^{b^1_P}_{Q'} = f \chi_{I_{Q_r}} - \frac{[f]_{Q_r}}{[b^1_P]_{Q_r}} b^1_P \chi_{I_{Q_r}}.$$
Subtracting the two we thus see that
$$ W_{b^1_P} f(Q) \psi^{b^1_P}_Q + \sum_{Q' \in T'} W_{b^1_P} f(Q') \psi^{b^1_P}_{Q'} = f \chi_{I_{T'}}
- \frac{[f]_Q}{[b^1_P]_{Q}} b^1_P \chi_{I_Q} + \frac{[f]_{Q_r}}{[b^1_P]_{Q_r}} b^1_P \chi_{I_{Q_r}}.$$
If we thus define
$$ \varphi_Q := [f]_{P_{T'}} b^1_{P_{T'}}
- \frac{[f]_Q}{[b^1_P]_{Q}} b^1_P \chi_{I_Q} + \frac{[f]_{Q_r}}{[b^1_P]_{Q_r}} b^1_P \chi_{I_{Q_r}}$$
we see that \eqref{aaaa} follows; the mean zero condition can be seen by \eqref{normal} and inspection.  This completes the proof of Lemma \ref{subtree}.
\end{proof}

Roughly speaking, the above lemma states that we can find a large tree $T_1$ on which $b^1_P$ is pseudo-accretive, on which $b^1_P$ and $\CZ b^1_P$ are effectively bounded, and for which we have a representation of the form \eqref{g-expand}.  We now run an argument in the spirit of Lemma \ref{tree-samp} to localize matters exclusively to this tree $T_1$.

We apply the above Lemma first with the $b^1_P$.  We decompose $f$ using \eqref{f-form}, thus estimating the left-hand side of \eqref{dual-targ} by the sum of the term below $T_1$
\be{dual-1}
|\int_{I_P} [f]_P \CZ b^1_P|,
\end{equation}
the terms coming from $T_1$
\be{dual-2}
|\sum_{Q \in T_1} W_{b^1_P} f(Q) \int_{I_P} \CZ \psi^{b^1_P}_Q |,
\end{equation}
the terms above $T_1 \cup \P_{buffer}$
\be{dual-3}
\sum_{T' \in \T} 
|\int_{I_P} \CZ(f \chi_{I_{T'}} - [f]_{P_{T'}} b^1_{P_{T'}})|,
\end{equation}
and the term from $\P_{buffer}$
\be{dual-4}
\sum_{Q \in \P_{buffer}} 
|\int_{I_P} \CZ \varphi_Q|.
\end{equation}

The contribution of \eqref{dual-1} is $O(|I_P| \|f\|_\infty)$ by \eqref{bounds} and H\"older.  For the contribution of \eqref{dual-2}, we observe from \eqref{f-ortho} that
$$ (\sum_{Q \in T_1} |W_{b^1_P} f(Q)|^2)^{1/2}
\lesssim \|f\|_2 \leq \|f\|_\infty |I_P|^{1/2}.$$
By Cauchy-Schwarz it will thus suffice to show the bound
$$
\sum_{Q \in T_1} |\int_{I_P} \CZ \psi^{b^1_P}_Q|^2 \lesssim |I_P|$$
or equivalently that
\be{tcarl}
\| |\langle \CZ^* \chi_{I_P}, \psi^{b^1_P}_Q \rangle|^2 \|_{\size(T_1)} \lesssim 1.
\end{equation}
One can think of \eqref{tcarl} as a localized, $b^1_P$-adapted version of the statement $\CZ^*(1) \in BMO$.

We will defer the proof of \eqref{tcarl} to Step 3 of the argument.  Assuming the bound \eqref{tcarl} for now, we move on to \eqref{dual-3}.  Observe that the expression inside the $\CZ()$ in \eqref{dual-3} is supported on $I_{T'}$ and has mean zero, so we may reduce the integral from $I_P$ to $I_{T'}$.  From the definition of $A$ we have
$$ |\int_{I_{T'}} \CZ(f \chi_{I_{T'}})| \leq A |I_{T'}| \|f\|_\infty$$
while from \eqref{bounds} and Cauchy-Schwarz we have
$$ |\int_{I_{T'}} \CZ([f]_{P_{T'}} b^1_{P_{T'}})| \lesssim |I_{T'}| \|f\|_\infty.$$
Adding all this up and using the fact that the tops of $\T$ are a $(1-\eps)$ packing we can bound \eqref{dual-3} by $((1-\eps)A + O(1))|I_P|$ as desired.

Finally we consider \eqref{dual-4}.  Since $\varphi_Q \in S_0(I_Q)$ we can estimate this term by
$$ \sum_{Q \in \P_{buffer}} |\int_{I_Q} \CZ \varphi_Q|.$$
If $|I_Q| = 2^{-M}$ then this vanishes because we truncated the operator $\CZ$, so we assume that $|I_Q| > 2^{-M}$.
Using Cauchy-Schwarz, \eqref{bounds}, \eqref{aaaa}, and Corollary \ref{trunc} we can bound this by
$$ O(\sum_{Q \in \P_{buffer}} |I_Q| \|f\|_\infty)$$
which is acceptable since $\P_{buffer}$ is a 2-packing.

{\bf Step 3.  Pruning the bad tiles of $b^2_P$.}

In Step 2 we reduced the proof of \eqref{dual-targ} (and thus of Theorem \ref{ctb} to that of proving the Carleson measure estimate \eqref{tcarl}.  Along the way we managed to prune all the tiles for which $b^1_P$ was ``bad'' (in that the mean of $b^1_P$ was too small, or the $L^2$ norm of $b^1_P$ or $\CZ b^1_P$ was too large).  However, $b^2_P$ is still not under control.  Thus the next step shall be to prune $b^2_P$.

We define
$$ B:=\| |\langle \CZ^* \chi_{I_P}, \psi^{b^1_P}_Q \rangle|^2 \|_{\size^*(T_1)};$$
it thus suffices to prove that $B = O(1)$.  In fact will suffice to show that
$$ \sum_{Q \in T_1 \cap \Tree(P')} |\langle \CZ^* \chi_{I_P}, \psi^{b^1_P}_Q \rangle|^2
\leq ((1 - \eps) B + O(1)) |I_{P'}|$$
for all tiles $P' \leq' P$, since the claim then follows by taking suprema over $P'$ and solving for $B$ (cf. Lemma \ref{tree-samp}).

Fix $P'$.  
We apply Lemma \ref{subtree} again but with the $b^2_{P'}$, partitioning
$$ \Tree(P') = T_2 \cup \P'_{buffer} \cup \bigcup_{T' \in \T'} T'.$$
From the definition of $B$ and the fact that $\T'$ is a $(1-\eps)$-packing we have (cf. Lemma \ref{tree-samp})
$$ \sum_{Q \in T_1 \cap \bigcup_{T' \in \T'} T'} |\langle \CZ^* \chi_{I_P}, \psi^{b^1_P}_Q\rangle|^2 \leq (1 - \eps) B |I_{P'}|,$$
so it suffices to show that
$$ \sum_{Q \in T_1 \cap (T_2 \cup \P'_{buffer})} |\langle \CZ^* \chi_{I_P}, \psi^{b^1_P}_Q\rangle |^2 \lesssim |I_{P'}|.$$
From \eqref{norm} it will suffice to show
$$
\sum_{Q \in T_1 \cap (T_2 \cup \P'_{buffer})} |\langle \CZ^* \chi_{I_P}, b^1_P \phi^{b^1_P}_Q\rangle |^2
\lesssim |I_{P'}|.
$$
We first observe that
$$ |\langle \CZ^* \chi_{I_P}, b^1_P \phi^{b^1_P}_Q \rangle|
= |\int_{I_P} \CZ (b^1_P \phi^{b^1_P}_Q)| =
|\int_{I_Q} \CZ (b^1_P \phi^{b^1_P}_Q)|
\leq |I_Q|^{1/2} \| \CZ (b^1_P \phi^{b^1_P}_Q) \|_2;$$
from Corollary \ref{trunc} and \eqref{mean-bounds} we thus have the weak Carleson bound
$$ |\langle \CZ^* \chi_{I_P}, b^1_P \phi^{b^1_P}_Q \rangle|^2 \lesssim |I_Q|$$
for all $Q \in T_1$.  In particular we see that the contribution of $\P'_{buffer}$ will be acceptable since $\P'_{buffer}$ is a $2$-packing.  We are thus left with showing
$$
\sum_{Q \in T_1 \cap T_2} |\langle \CZ^* \chi_{I_P}, b^1_P \phi^{b^1_P}_Q\rangle|^2
\lesssim |I_{P'}|.
$$
In Step 4 we shall show the pointwise bound
\be{pointwise}
|\langle \CZ^* \chi_{I_P}, b^1_P \phi^{b^1_P}_Q\rangle|
\lesssim 
|Wb^2_{P'}(Q)| + |W_{b^1_P}(b^1_P \CZ^*(b^2_{P'}))(Q)|
+ |W_{b^1_P}(b^2_{P'} \CZ(b^1_P))(Q)|
+ |W_{b^1_P}(\CZ(b^1_P))(Q)|
\end{equation} 
for all $Q \in T_1 \cap T_2$.  The claim will then follow from \eqref{f-ortho}, \eqref{f-var}\footnote{While the tree $T_1 \cap T_2$ is not necessarily convex, the larger tree $(T_1 \cup \P_{buffer}) \cap (T_2 \cup \P'_{buffer})$ is, and on this larger tree we have pseudo-accretivity and \eqref{mean-bounds} for both $b^1_P$ and $b^2_{P'}$.}.

{\bf Step 4.  Pointwise estimates on wavelet coefficients.}

In the previous step we reduced matters to proving the pointwise estimate \eqref{pointwise}.  The proof of this estimate is really the core of our argument, although it was necessary to do all the above prunings to get to a point where this estimate became both provable and useful\footnote{If one wished to prove a global $\CZ(b)$ theorem, with globally para-accretive $L^\infty$ functions $b^1_P$, $b^2_P$, one could dispense with the selection algorithms and go directly to (a suitable analogue of) \eqref{pointwise}.  We omit the details.}.

We now prove \eqref{pointwise}.  Fix $Q \in T_1 \cap T_2$.  On $I_Q$, we can decompose
$ \chi_{I_P} = \frac{b^2_{P'}}{[b^2_{P'}]_Q} + F$
where $F$ is the mean zero function
$F := \chi_{I_Q} - \frac{b^2_{P'} \chi_{I_Q}}{[b^2_{P'}]_Q}$.
Since $\CZ(b^1_P \phi^{b^1_P}_Q)$ is supported on $I_Q$, we can thus estimate the left-hand side of \eqref{pointwise} by
$$ |\langle \CZ^* \frac{b^2_{P'}}{[b^2_{P'}]_Q}, b^1_P \phi^{b^1_P}_Q \rangle|
+ |\langle \CZ^* F, b^1_P \phi^{b^1_P}_Q \rangle|.$$
The first term is $|O(W_{b^1_P}(b^1_P \CZ^*(b^2_{P'}))(Q))|$ by the pseudo-accretivity of $b^2_{P'}$.  For the second term, observe that if the $\phi^{b^1_P}_Q$ could be moved inside the $\CZ^*$, thus
$$ |\langle \CZ^* (\phi^{b^1_P}_Q F), b^1_P \rangle|$$
then by moving the $\CZ^*$ to the other side, we could bound this by
$$ O(|W_{b^1_P}(F \CZ(b^1_P))|) = O(|W_{b^1_P}(b^2_{P'} \CZ(b^1_P))(Q)|) + O(|W_{b^1_P}(\CZ(b^1_P))(Q)|)$$
again using the pseudo-accretivity of $b^2_{P'}$.  Thus it suffices to control the commutator
$$ | \langle \CZ^* F, b^1_P \phi^{b^1_P}_Q \rangle - \langle \CZ^* (\phi^{b^1_P}_Q F), b^1_P \rangle|.$$
If $F$ had mean zero on both $I_{Q_l}$ and $I_{Q_r}$ then this commutator would be zero from \eqref{perfect} since $\phi^{b^1_P}_Q$ is constant on $I_{Q_l}$ and $I_{Q_r}$.  Thus we may freely replace $F$ by $[F]_{Q_l} (\chi_{I_{Q_l}} - \chi_{I_{Q_r}})$ since the difference has mean zero on both $I_{Q_l}$ and $I_{Q_r}$.  Throwing the $\CZ^*$ on to the other side and using Cauchy-Schwarz and Corollary \ref{trunc}, we can thus bound this commutator by
$$ O( |I_Q|^{1/2} |[F]_{Q_l}|).$$
However a computation shows that
$$ [F]_{Q_l} = \frac{1}{2} ([F]_{Q_l} - [F]_{Q_r}) =  \frac{[b^2_{P'}]_{Q_r} - [b^2_{P'}]_{Q_l}}{2[b^2_{P'}]_Q} = -|I_Q|^{-1/2} \frac{Wb^2_{P'}(Q)}{[b^2_{P'}]_Q}$$
and the claim then follows by the pseudo-accretivity of $b^2_{P'}$.  The proof of Theorem \ref{ctb} is now complete.
\end{proof}

One can generalize \eqref{bounds} to
$$
\int_{I_P} |b^1_P|^p + |\CZ b^1_P|^{q'} + |b^2_P|^q + |\CZ^*b^2_P|^{p'} \lesssim |I_P|$$
for any $1 \leq p,q \leq \infty$, with $1/p + 1/p' = 1/q + 1/q' = 1$ (this was already suspected in \cite{christ-tb} when $p=q=\infty$); the dual exponents are necessary to control such expressions as $\langle \CZ b^1_P, b^2_P \rangle$.  Most of the argument proceeds \emph{mutatis mutandis} except for Lemma \ref{ortho-lemma}.  Firstly in \eqref{b-bmo} the $|b|^2$ mean has to be replaced by some other $|b|^p$ mean, but because of Corollary \ref{jn-char} we still recover \eqref{b-bmo-2}.  However we still must modify \eqref{f-ortho}, \eqref{f-var} to
$$
(\sum_{P \in T} |W_b f(P)|^2)^{1/2} \lesssim |I_T|^{1/2} \| |f|^{q'} \|_{\mean^*(T)}^{1/q'}
$$
and
$$
(\sum_{P \in T} |W_b (b' f) (P)|^2)^{1/2} \lesssim |I_T|^{1/2} \| |f|^{q'} \|_{\mean^*(T)}^{1/q'} \| |b'|^q \|_{\mean^*(T)}^{1/q}.$$
This proceeds by replacing $f$ and $b' f$ with averaged variants as in the proof of \eqref{f-var}; the averages will then be controlled in $L^\infty$ and hence in $L^2$.  We omit the details.

It is also straightforward to generalize Theorem \ref{ctb} to Calder\'on-Zygmund operators which do not obey the perfect dyadic cancellation condition \eqref{perfect}, and instead obey a more classical cancellation condition such as $|\nabla_x K(x,y)| + |\nabla_y K(x,y)| \lesssim 1/|x-y|^2$.  However it is still convenient to impose a truncation condition on the kernel when $|x-y|$ is extremely small or extremely large (as in e.g. \cite{christ-tb}).  

The main new difference with these kernels is that when $f_I \in S_0(I)$, the function $\CZ f_I$ is no longer supported in $I$ but has a tail at infinity.  However the cancellation conditions ensure that this tail is quite rapidly decaying (like $1/|x|^2$ if we assume the above gradient bounds).  This causes many of the identities used in the arguments above to pick up some error terms, for instance if $g$ is constant on $I$ it is no longer true that $\langle \CZ f_I, g \rangle = [g]_I \langle \CZ f_I, 1 \rangle$, however the error term incurred is quite manageable due to the good decay (especially if $g$ is in fact constant on a much wider interval than $I$).  We will not pursue the details further here as they are rather standard (see e.g. \cite{cjs}, \cite{christ-tb}).  A perhaps more interesting generalization would be to non-doubling situations as in \cite{ntv:tb-local}, as this may have applications to analytic capacity problems, but this seems to require much more technical arguments.

As a corollary of the local $\CZ(b)$ theorem we can conclude a global $\CZ(b)$ theorem\footnote{Of course, the global $\CZ(b)$ theorem could be proven directly in a much simpler manner, but one advantage of doing things this way is that we can relax the hypotheses of the global $\CZ(b)$ theorem slightly.  If $b_1$, $b_2$ were bounded and strongly pseudo-accretive, one could obtain a direct proof of the global $\CZ(b)$ theorem by using the modified wavelet transforms $W_{b_1}$, $W_{b_2}$ to define paraproducts by adapting \eqref{permute}, and then finding an analogue of \eqref{decomp}, and repeating the proof of the global $\CZ(1)$ theorem.  See e.g. \cite{at} for details.  An alternate approach based on \eqref{pointwise} is also possible.} .  We recall that a function $b$ is \emph{para-accretive} if for every tile $P$ there exists a tile $Q \leq' P$ with $|I_Q| \sim |I_P|$ and $|[b]_Q| \gtrsim 1$.  Every pseudo-accretive function
is para-accretive (just take $Q := P$), but not conversely.

\begin{corollary}[Dyadic global $\CZ(b)$ theorem]\label{tb}
Let $b_1, b_2$ be para-accretive functions with 
$$\| b_1\|_{BMO}, \|b_2\|_{BMO} \lesssim 1.$$
Let $\CZ$ be a perfect dyadic Calder\'on-Zygmund operator obeying the modified weak boundedness property
\be{mwbp}
|\langle \CZ b_1 \chi_I, b_2 \chi_J \rangle| \lesssim 
|K| \hbox{ for all intervals } I, J, K \hbox{ with } |I| \sim |J| \sim |K| \hbox{ and } I, J \subseteq K
\end{equation}
and such that 
$$ \| \CZ(b_1)\|_{BMO}, \|\CZ^*(b_2) \|_{BMO} \lesssim 1.$$
Then $\CZ$ is bounded on $L^2$.
\end{corollary}

\begin{proof}
We apply Theorem \ref{ctb} with
$$ b^1_P := \frac{b_1 \chi_{I_Q}}{[b_1]_Q}$$
where for each $P$, we choose the tile $Q \leq' P$ so that $|I_Q| \sim |I_P|$ and $|[b_1]_Q| \gtrsim 1$.  We define $b^2_P$ similarly.

The normalization \eqref{normal} is clear.  To prove \eqref{bounds}, we observe that
$$ b^1_P = \frac{(b_1 - [b_1]_Q) \chi_{I_Q}}{[b_1]_Q} + 1$$
and so the $L^2$ bound on $b^1_P$ follows from the BMO control on $b_1$ and the lower bound on $|[b_1]_Q|$.  To control $\CZ b^1_P$, it suffices from \eqref{kernel} and the lower bound on $|[b_1]_Q|$ to show that
$$ \| \CZ (b_1 \chi_{I_Q}) \|_{L^2(I_Q)} \lesssim  |I_Q|^{1/2},$$
or in other words that
$$ |\langle \CZ(b_1 \chi_{I_Q}), h \rangle| \lesssim |I_Q|^{1/2} \|h\|_2$$
for all $h \in S(I_Q)$.

Select a tile $R \leq' Q$ such that $|[b_2]_R| \gtrsim 1$ and $|I_R| \sim |I_Q|$.  If $h$ is a scalar multiple of $b_2 \chi_R$, the claim follows from \eqref{mwbp}.  Thus we may subtract multiples of $b_2 \chi_R$ and reduce to the case when $h$ has mean zero (cf. Lemma \ref{split-lemma}).  But then $\langle \CZ(b_1 \chi_{I_Q}), h \rangle = \langle \CZ b_1, h \rangle$, and the claim follows since $\CZ b_1 \in BMO$. 

One can control $b^2_P$ and $\CZ^* b^2_P$ by identical arguments.
\end{proof}

Notice that $b_1$, $b_2$ are only assumed to be in BMO\footnote{For untruncated Calder\'on-Zygmund operators in the limit $M \to \infty$ it may well be necessary to require $b_1$, $b_2$ to be in VMO rather than BMO.} rather than $L^\infty$.  This generalization of the standard $\CZ(b)$ theorem appears to be new.  Also observe that the above argument also works in the special case $\CZ(b_1) = \CZ^*(b_2)=0$ if we drop the para-accretivity and BMO hypotheses on $b_1$, $b_2$ and instead impose the reverse H\"older conditions
$$ \frac{(\frac{1}{|I_P|} \int_{I_P} |b_1|^2 )^{1/2}}{|[b_1]_P|},
\frac{(\frac{1}{|I_P|} \int_{I_P} |b_2|^2 )^{1/2}}{|[b_2]_P|}
\lesssim 1$$
on $b_1$, $b_2$.  It is in fact likely that we can obtain a $\CZ(b)$-type theorem for arbitrary (complex) dyadic $A_\infty$ weights $b_1$, $b_2$ (see \cite{fkp}), but we will not attempt to give the most general statements here.

\end{document}